\documentclass[12pt]{article}

\usepackage{amssymb}
\usepackage{amsthm}
\usepackage{amsmath}
\usepackage[T1]{fontenc}
\usepackage[utf8]{inputenc}
\usepackage[english]{babel}
\usepackage[pdftex]{graphicx}      
\usepackage{mathrsfs}

\usepackage{indentfirst}              
\usepackage{makeidx}                    
\usepackage[nottoc]{tocbibind}         
\usepackage{courier}                   
\usepackage{type1cm}                    
\usepackage{listings}                
\usepackage{titletoc}
\usepackage[all,cmtip]{xy}

\usepackage{fancyhdr}
\newtheorem{defi}{\bf Definition}[section]
\newtheorem{theo}[defi]{\bf Theorem}

\newtheorem{coro}[defi]{\bf Corollary}
\newtheorem{pro}[defi]{\bf Proposition}
\newtheorem{lem}[defi]{\bf Lemma}

\newenvironment{dem}{\noindent
\textnormal{\textbf{Proof:}}}{\begin{flushright}$\mathbf{Q.E.D.}$\end{flushright}}

\theoremstyle{remark}

\newcommand{\algA}{\mathcal{A}}

\title{Kronecker factorization theorems for the exceptional Malcev algebra}

\author{
{Victor L\'opez Sol\'{\i}s
}\\
{\small Departamento Académico de Ciencias Básicas y Afines}\\
{\small Universidad Nacional de Barranca,}\\
{\small Jr. Toribio de Luzuriaga N°\,376, Barranca, 15169},\\
{\small Lima. Peru}\\
{\small vlopez@unab.edu.pe}
}
\date{\quad}

\begin{document}
\maketitle\vspace{-1.5cm}


\begin{abstract}

We prove that a Malcev algebra $\mathcal{M}$ which belongs to a certain variety of Malcev algebras $\mathcal{H}$ such that $\mathcal{M}$ contains the $7$-dimensional exceptional Malcev algebra $\mathbb{M}$ is isomorphic to a direct sum $(\mathbb{M}\otimes_{F} A)\oplus \mathcal{N}$, where $A$ is a commutative associative algebra.
Also, we prove that a Malcev superalgebra $\mathcal{M}=\mathcal{M}_0\oplus \mathcal{M}_1$ in $\mathcal{H}$ whose even part $\mathcal{M}_0$ contains $\mathbb{M}$ is isomorphic to a direct sum $(\mathbb{M}\otimes_{F} \mathcal{A})\oplus \mathcal{N'}$, where $\mathcal{A}$ is a  supercommutative associative superalgebra.
\end{abstract}

{\parindent= 4em \small  \sl Keywords: Malcev (super)algebra,  non-Lie Malcev algebra, Malcev module, Kronecker factorization theorem.\\
MSC classification (2020): Primary 17D10, secondary 17B22, 17B60, 17B70.\\
https://orcid.org/
0000-0001-5969-5853 
}


\section{Introduction}

A \textit{Malcev algebra} is a vector space $\mathcal{M}$ with a bilinear binary operation $(x,y)\mapsto xy $ satisfying the following identities:
\begin{equation}\label{e}
x^2=0, ~~ J(x,y,xz)=J(x,y,z)x, 
\end{equation}
where $J(x,y,z)=(xy)z+(yz)x+(zx)y$ is the Jacobian of the elements $x,y,z\in \mathcal{M}$. Denote by $J(\mathcal{M},\mathcal{M},\mathcal{M})$ the subspace spanned by the Jacobians, which is an ideal of $\mathcal{M}$ (see $\cite{SA}$). The set $N_{\mathcal{M}}=\{x \in \mathcal{M}: J(x, \mathcal{M}, \mathcal{M})=0\}$ is also an ideal of $\mathcal{M}$ called the \textit{Lie nucleus} of $\mathcal{M}$. A Malcev algebra $\mathcal{M}$ is a Lie algebra if and only if $J(\mathcal{M}, \mathcal{M}, \mathcal{M})=0$ if and only if $N_{\mathcal{M}}=\mathcal{M} .$ In some way, the Malcev algebras with trivial Lie nucleus, $N_{\mathcal{M}}=0$, are the most distant from being Lie algebras.

All Lie algebras are clearly Malcev algebras because the Jacobian of any three elements vanish. The tangent space $T(L)$ of an analytic Moufang loop $L$ is another example of Malcev algebra. Let $\algA$ be an alternative algebra, if we introduce a new product by means of a commutator $[x,y]=xy-yx$ into $\algA$, we obtain a new algebra that will be denoted by $\algA^{(-)}$. It is easy to verify that the algebra $\algA^{(-)}$ is a Malcev algebra. All Malcev algebras obtained in this form are called \textit{special}. 
The classic example of a non-Lie Malcev algebra $\mathbb{M}$ is formed by traceless elements of the Cayley-Dickson algebra with the commutator. The algebra $\mathbb{M}$ is of dimension $7$ and is one of most cornerstone examples.

H. M. Wedderburn proved that if a unital associative algebra $\algA$ contains a central simple subalgebra of finite dimension $\mathcal{B}$ with the same identity element, then $\algA$ is isomorphic to a Kronecker product $S\otimes_F \mathcal{B}$, where $S$ is the subalgebra of the elements that commute with each $b\in \mathcal{B}$. In particular, if $\algA$ contains $M_{n}(F)$ as subalgebra with the same identity element, we have $\algA\cong S\otimes_F M_{n}(F)$ where $S$ is the subalgebra of the elements that commute with each unitary matrix $e_{ij}$. 
Kaplansky extended (see \cite{K}, Theorem $2$) Wedderburn result to alternative algebras and the split Cayley algebra.
Jacobson uses his classification of completely reducible alternative bimodules over fields of characteristic diferent of 2 to give a new proof of Kaplansky result in \cite{J1}, Theorem $1$. The author proved that the result is valid for any characteristic (see \cite{LSS}).
Jacobson used his above mentioned result in \cite{J1} to prove a Kronecker Factorization Theorem for the simple exceptional Jordan algebra of $3\times 3$ Hermitian matrices over the Cayley algebra.
In the case of right alternative algebras, S. Pchelintsev, O. Shashkov and I. Shestakov \cite{PSS} proved that every unital right alternative bimodule over a Cayley algebra (over an algebraically closed field of characteristic not 2) is alternative and they used that result to prove a coordinatization theorem for unital right alternative algebras containing a Cayley subalgebra with the same unit. 

In the case of superalgebras, M. López-Díaz and I. Shestakov \cite{LDS2, LDS1} studied the representations of simple alternative and exceptional Jordan superalgebras in characteristic $3$ and through these representations, they obtained some analogues of the Kronecker Factorization Theorem for these superalgebras. 
Also, the author \cite{LSS} obtained some analogues of the Kronecker Factorization Theorem for central simple alternative superalgebras $M_{(1\mid 1)}(F)$,  $\mathbb{O}(4|4)$ and $\mathbb{O}[u]$, where in particular the Kronecker Factorization Theorem for $M_{(1\mid 1)}(F)$ answers the analogue for superalgebras of the Jacobson's problem \cite{J1}. We state that Jacobson's problem \cite{J1} was recently solved by the author and I. Shestakov \cite{LSS1} in the split case.
Similarly, C. Martinez and E. Zelmanov \cite{CE} obtained a Kronecker Factorization Theorem for the exceptional ten dimensional Kac superalgebra $K_{10}$. Also, A. Pozhidaev and I. Shestakov  \cite{PS} studied the representations of simple finite-dimensional noncommutative Jordan superalgebras and proved some analogues of the Kronecker factorization theorem for such superalgebras (see also \cite{PY}). 

In \cite[(6)]{F2}, Filippov introduced a very interesting identity and considered (see \cite{F4,F,F1,F3}) the subvariety $\mathcal{H}$ of those Malcev algebras satisfying this identity. He proved in \cite{F} that, in a sense, the subvariety $\mathcal{H}$ plays the same role as the variety of Lie algebras in the structure theory of Malcev algebras. Moreover, the central closures of the prime algebras in $\mathcal{H}$ and the simple algebras in $\mathcal{H}$ are the 3-dimensional simple Lie algebra $\mathfrak{sl}(2,F)$ and the 7-dimensional exceptional Malcev algebra $\mathbb{M}$ (see \cite{ES1}).

The aim of this paper is to determine the structure of Malcev (super)algebras in $\mathcal{H}$ containing the $7$-dimensional exceptional Malcev algebra $\mathbb{M}$. So the paper is organized as follows: In Section 2, we provided definitions about Malcev (super)algebras and their representations. 
Also we describe the special subvariety $\mathcal{H}$ of Malcev (super)algebras and their representations. 
In Section 3, we prove some Kronecker Factorization Theorems for Malcev (super)algebras in $\mathcal{H}$ that contain $\mathbb{M}$. In Section
4, we establish certain equivalences of categories.
Finally, in section $5$, we study (super)algebras and modules with (super)involutions.


\section{Preliminaries} 

In this section we provide background material that is used along the way and some preliminary results.

For arbitrary elements $x,y,z,t,u$ in a Malcev algebra $\mathcal{M}$ we consider the following functions that play an important role in
the theory of Malcev algebras:
$$J(x,y,z)=(xy)z+(yz)x+(zx)y,~\mbox{\textit{the Jacobian} of}~x,y,z,$$
$$[x,y,z]=(xy)z+x(yz),~\mbox{\textit{the antiassociator}},$$
$$\{x,y,z\}=(xy)z-(xz)y+2x(yz)=J(x,y,z)+3x(yz)=[x,y,z]-[x,z,y],$$
$$h(y,z,t,x,u)=\{yz,t,u\}x+\{yz,t,x\}u+\{yx,z,u\}t+\{yu,z,x\}t.$$

Throughout this paper $F$ will be a field of characteristic different of 2 and 3. Then expanding the Jacobian, the Malcev identities $(\ref{e})$ can be rewritten as follows
\begin{equation}\label{line1}
\begin{split}
xy&=-yx,\\
(xz)(yt)=((xy)z)t+((&yz)t)x+((zt)x)y+((tx)y)z.
\end{split}
\end{equation}
And the following identity hold in $\mathcal{M}$ which was proved in $\cite{SA}$
\begin{equation}\label{e8}
    J(tx,y,z)=tJ(x,y,z)+J(t,y,z)x-2J(t,x,yz).
\end{equation}

\subsection{Seven dimensional exceptional Malcev algebra}

The most important example of a non-Lie Malcev algebra is the algebra $\mathbb{M}$ of dimension $7$ which we are going to describe now. First, we consider the central simple Malcev algebra $\mathbb{M}=\mathcal{M}_7(F)$ 
obtained by means of the commutator $[x,y]=xy-yx$ introduced in the \textit{split Cayley-Dickson matrix algebra} $C(F)$. 

As in $\cite{SA}$, we write the elements of the algebra $C(F)$ in the form
$$
\begin{bmatrix}
a&\alpha\\
\beta&b\\
\end{bmatrix},
$$
where $a,b$ are elements of the field $F$ and $\alpha,\beta$ are $3$-dimensional vectors $(a_{1},a_{2},a_{3})$ over $F$.

The sum is defined as
$$
\begin{bmatrix}
a&\alpha\\
\beta&b\\
\end{bmatrix}
+
\begin{bmatrix}
c&\gamma\\
\delta& d\\
\end{bmatrix}
=
\begin{bmatrix}
a+c&\alpha+\gamma\\
\beta+\delta&b+d\\
\end{bmatrix}.
$$

We use the notations $\bullet$ and $\times$ for the scalar and product of $3$-dimensional vectors and define the multiplication by
$$
\begin{bmatrix}
a&\alpha\\
\beta&b\\
\end{bmatrix}
.
\begin{bmatrix}
c&\gamma\\
\delta& d\\
\end{bmatrix}
=
\begin{bmatrix}
ac+\alpha\bullet\delta&a\gamma+d\alpha-\beta\times\delta\\
c\beta+b\delta+\alpha\times\gamma&\beta\bullet\gamma+bd\\
\end{bmatrix}.
$$

In the split Cayley-Dickson algebra $C(F)$ we define a new multiplication $[\,,\,]$ by $[A,B]=AB-BA$, or
\begin{eqnarray*}
[A,B]&=&
\begin{bmatrix}
a&\alpha\\
\beta&b\\
\end{bmatrix}
.
\begin{bmatrix}
c&\gamma\\
\delta& d\\
\end{bmatrix}
-
\begin{bmatrix}
c&\gamma\\
\delta& d\\
\end{bmatrix}
.
\begin{bmatrix}
a&\alpha\\
\beta&b\\
\end{bmatrix}\\
&=&
\begin{bmatrix}
\alpha\bullet\delta-\beta\bullet\gamma&(a-b)\gamma+(d-c)\alpha-2\beta\times\delta\\
(c-d)\beta+(b-a)\delta+2\alpha\times\gamma&-(\alpha\bullet\delta-\beta\bullet\gamma)\\
\end{bmatrix}.
\end{eqnarray*}

The algebra $C(F)^{(-)}$, with multiplication $[\,,\,]$ defined above, is a Malcev algebra. We chose a basis for $C(F)^{(-)}$ the elements
$$
e_{0}=
\begin{bmatrix}
1&0\\
0&1\\
\end{bmatrix},
~~~e_{1}=
\begin{bmatrix}
1&0\\
0&-1\\
\end{bmatrix},
~~~e_{2}=
\begin{bmatrix}
0&(1,0,0)\\
0&0\\
\end{bmatrix},
$$
$$
e_{3}=
\begin{bmatrix}
0&(0,1,0)\\
0&0\\
\end{bmatrix},
~~~e_{4}=
\begin{bmatrix}
0&(0,0,1)\\
0&0\\
\end{bmatrix},
~~~e_{5}=
\begin{bmatrix}
0&0\\
(1,0,0)&0\\
\end{bmatrix},
$$
$$
e_{6}=
\begin{bmatrix}
0&0\\
(0,1,0)&0\\
\end{bmatrix},
~~~e_{7}=
\begin{bmatrix}
0&0\\
(0,0,1)&0\\
\end{bmatrix}.
$$

For all $x$ in $C(F)^{(-)}$, we have $[e_{0},x]=0$ and so let us consider the algebra 
$$\mathcal{M}_7 (F)=C(F)^{(-)}/F\cdot e_{0}.$$ 
We identify $e_{i}$ in $C(F)^{(-)}$ with the coset $e_{i}+F\cdot e_{0}$ in $\mathcal{M}_7 (F)$ and we denote the multiplication $[\,,\,]$ by juxtaposition, that is, $xy=[x,y]$. So we have the following multiplication table of $\mathcal{M}_7 (F)$ 

\begin{center}
\begin{tabular}{|c||c|c|c|c|c|c|c|}
\hline
$.$ & $e_{1}$ & $e_{2}$ & $e_{3}$ & $e_{4}$ & $e_{5}$ & $e_{6}$ & $e_{7}$\\\hline\hline
$e_{1}$ & $0$ & $2e_{2}$ & $2e_{3}$ & $2e_{4}$ & $-2e_{5}$ & $-2e_{6}$ & $-2e_{7}$\\\hline
$e_{2}$ & $-2e_{2}$ & $0$ & $2e_{7}$ & $-2e_{6}$ & $e_{1}$ & $0$ & $0$\\\hline
$e_{3}$ & $-2e_{3}$ & $-2e_{7}$ & $0$ & $2e_{5}$ & $0$  & $e_{1}$ & $0$\\\hline
$e_{4}$ & $-2e_{4}$ & $2e_{6}$ & $-2e_{5}$ & $0$ & $0$ & $0$ & $e_{1}$\\\hline
$e_{5}$ & $2e_{5}$ & $-e_{1}$ & $0$ & $0$ & $0$ & $-2e_{4}$ & $2e_{3}$\\\hline
$e_{6}$ & $2e_{6}$ & $0$ & $-e_{1}$ & $0$ & $2e_{4}$ & $0$ & $-2e_{2}$\\\hline 
$e_{7}$ & $2e_{7}$ & $0$ & $0$ & $-e_{1}$ & $-2e_{3}$ & $2e_{2}$ & $0$\\\hline
\end{tabular}
\end{center}
and we know that the algebra $\mathcal{M}_7 (F)$ is simple.

Now we consider the Malcev algebra $\mathbb{M}=\mathbb{M}_7$ obtained by means of the commutator $[x,y]=xy-yx$ introduced in \textit{the division Cayley-Dickson algebra} $\mathbb{O}$. Remember that $\mathbb{O}=\mathbb{H}\oplus v\mathbb{H}$, $v^{2}=\gamma\neq 0$ ($\gamma\in F$) has the $\mathbb{Z}_{2}$-grading induced by the Cayley-Dickson process applied to the division quaternion algebra $\mathbb{H}$. So, the odd part $v\mathbb{H}$ is a bimodule over $\mathbb{H}$ and it is a \textit{Cayley bimodule} because satisfy the relation
\begin{equation}\label{cay1}
aw=w\bar{a},
\end{equation}
where $a\in\mathbb{H}$, $w\in v\mathbb{H}$ and $a\rightarrow \bar{a}$ is the canonical involution in $\mathbb{H}$. Also, we have
\begin{equation}\label{cay2}
\begin{split}
(ab)v&=b(av),\\
v(ab)&=(vb)a,\\
(va)(vb)&= \gamma b\bar{a},
\end{split}
\end{equation}
for all $a,b\in\mathbb{H}$. 

The algebra $\mathbb{O}^{(-)}$, with multiplication $[x,y]=xy-yx$, is a Malcev algebra. As the characteristic of the base field is different of $2$, let $\left\{1,i,j,k,v,vi,vj,vk\right\}$ be a basis of $\mathbb{O}^{(-)}$, where
$$i^{2}=j^{2}=k^{2}=-1.$$
For all $x$ in $\mathbb{O}^{(-)}$, we have $[1,x]=0$ and so let us consider the algebra 
$$\mathbb{M}_7=\mathbb{O}^{(-)}/F\cdot 1.$$
We identify any the base elements $x\in\left\{1,i,j,k,v,vi,vj,vk\right\}$ in $\mathbb{O}^{(-)}$ with the coset $x+F\cdot 1$ in $\mathbb{M}_7$ and again we denote the multiplication $[\,,\,]$ by juxtaposition, that is, $xy=[x,y]$. Then by $(\ref{cay1})$ and $(\ref{cay2})$ we have the following multiplication table of $\mathbb{M}_7$ 

\begin{center}
\begin{tabular}{|c||c|c|c|c|c|c|c|}
\hline
$.$ & $i$ & $j$ & $k$ & $v$ & $vi$ & $vj$ & $vk$\\\hline\hline
$i$ & $0$ & $2k$ & $-2j$ & $-2vi$ & $2v$ & $-2vk$ & $2vj$\\\hline
$j$ & $-2k$ & $0$ & $2i$ & $-2vj$ & $2vk$ & $2v$ & $-2vi$\\\hline
$k$ & $2j$ & $-2i$ & $0$ & $-2vk$ & $-2vj$ & $2vi$ & $2v$\\\hline
$v$ & $2vi$ & $2vj$ & $-2vk$ & $0$ & $2\gamma i$ & $2\gamma j$ & $2\gamma k$\\\hline
$vi$ & $-2v$ & $-2vk$ & $2vj$ & $-2\gamma i$ & $0$ & $2\gamma k$ & $-2\gamma j$\\\hline
$vj$ & $2vk$ & $-2v$ & $-2vi$ & $-2\gamma j$ & $-2\gamma k$ & $0$ & $2\gamma i$\\\hline 
$vk$ & $-2vj$ & $2vi$ & $-2v$ & $-2\gamma k$ & $2\gamma j$ & $-2\gamma i$ & $0$\\\hline
\end{tabular}
\end{center}
and also we know that the algebra $\mathbb{M}_7$ is simple.

\subsection{Representation of Malcev superalgebras}

In this subsection we provide some definitions about the representation theory of Malcev superalgebras. 

Let $\mathcal{M}=\mathcal{M}_0\oplus \mathcal{M}_1$ be a superalgebra (this is, $\mathbb{Z}_2$-graded algebra) and let $G=G_0\oplus G_1$ be the Grassmann algebra generated by the elements $1,e_{1},\ldots,e_{n},\ldots$ over a field $F$. The \textit{Grassmann envelope} of $\mathcal{M}$ is defined to be $G(\mathcal{M}):=G_0\otimes \mathcal{M}_0+G_1\otimes \mathcal{M}_1$. Then $\mathcal{M}$ is said to be \textit{Malcev superalgebra} if $G(\mathcal{M})$ is a Malcev algebra. So from this definition and by $(\ref{line1})$, $\mathcal{M}$ satisfies the following superidentities
\begin{equation*}
    xy=-(-1)^{|x||y|}yx,
\end{equation*}
\begin{equation*}
\begin{split}
(-1)^{|y||z|}(xz)(yt)&=((xy)z)t+(-1)^{|x|(|y|+|z|+|t|)}((yz)t)x\\
                    &\quad +(-1)^{(|x|+|y|)(|z|+|t|)}((zt)x)y\\
					&\quad +(-1)^{(|x|+|y|+|z|)|t|}((tx)y)z.
\end{split}
\end{equation*}
where $x,y,z\in\mathcal{M}_0\cup\mathcal{M}_1$ and $|x|$ denotes the parity index of a homogeneous
element $x$ of $\mathcal{M}$: $|x|=i$ if $x\in\mathcal{M}_i$.

In any Malcev superalgebra  $\mathcal{M}=\mathcal{M}_0\oplus \mathcal{M}_1$, the identities for Malcev algebras are easily superized to obtain the analogous graded identities. Thus, for homogeneous elements $x,y,z,t,u\in\mathcal{M}$ we define the following functions 
\begin{center}
    $\Tilde{J}(x,y,z)=(xy)z-x(yz)-(-1)^{|y||z|}(xz)y,\,\,\,\,\, (super\, Jacobian)$
\end{center}
\begin{equation*}
\begin{split}
\{x,y,z\}^{\tilde{}}&=(xy)z-(-1)^{|y||z|}(xz)y+2x(yz),\\
\Tilde{h}(y,z,t,x,u)&=\{yz,t,u\}^{\tilde{}}x+(-1)^{|x||u|}\{yz,t,x\}^{\tilde{}}u\\
&\,\,\,\,\,\,\,+(-1)^{|x|(|z|+|t|+|u|)+|u||t|}\{yx,z,u\}^{\tilde{}}t+(-1)^{|u|(|z|+|t|)+|x||t|}\{yu,z,x\}^{\tilde{}}t.
\end{split}
\end{equation*}

If $V=V_0\oplus V_1$ is a $\mathbb{Z}_2$-graded $F$-space, then $\mbox{End}_F(V)$ becames an associative superalgebra with
$$ \mbox{End}_F(V)_0=\{\varphi\in \mbox{End}_F(V):\varphi(V_i)\subseteq V_{i},\,i=0,1\},$$
$$ \mbox{End}_F(V)_1=\{\varphi\in \mbox{End}_F(V):\varphi(V_i)\subseteq V_{1-i},\,i=0,1\}.$$
The space $V$ is said to be \textit{module} for the Malcev superalgebra $\mathcal{M}$ if there is an even $F$-linear map $\rho:\mathcal{M}\longrightarrow \mbox{End}_F(V)(x\mapsto  \rho_{x})$, such that the split null extension $E=\mathcal{M}\oplus V$, with $\mathbb{Z}_2$-grading given by $E_0=\mathcal{M}_0\oplus V_0$, $E_1=\mathcal{M}_1\oplus V_1$, and with multiplication determined for homogeneus elements by
\begin{equation}\label{super1}
    (x+v)(y+w)=xy+(\rho_{y}(v)-(-1)^{|x||w|}\rho_{x}(w))
\end{equation}
is a Malcev superalgebra. The map $\rho$ is called a \textit{(super)-representation} of $\mathcal{M}$.

When dealing with superalgebras, the ideals or submodules will always be graded. The module $V$ is \textit{irreducible} if $\rho\neq 0$ and does not contain any proper submodule. Also, $V$ is said to be \textit{almost faithful} if $\mbox{ker}\rho$ does not contain any nonzero ideal of $\mathcal{M}$.

A \textit{regular module}, $\mbox{Reg}\,\mathcal{M}$, for superalgebra $\mathcal{M}$, is defined on the vector space $\mathcal{M}$ with the action of $\mathcal{M}$ coinciding with the multiplication in $\mathcal{M}$.

Let $V$ be a module for the Malcev superalgebra $\mathcal{M}$ and let $E=\mathcal{M}\oplus V$ be the corresponding split null extension. Let us consider
$$\widetilde{\Gamma}_i=\{\alpha\in\mbox{End}_F(E)_i:(xy)\alpha=x(y\alpha)=(-1)^{i|y|}(x\alpha)y\,\,\forall~x,y\in E_0\cup E_1\},\, i=0,1,$$
$$\widetilde{\Gamma}=\widetilde{\Gamma}(E)=\widetilde{\Gamma}_0\oplus \widetilde{\Gamma}_1,\,\,\mbox{\textit{the supercentroid} of}\,\,E,$$
$$Z=Z(E)=\{\alpha\in\widetilde{\Gamma}:~V\alpha\subseteq V~\mbox{and}~\mathcal{M}\alpha\subseteq \mathcal{M}\},$$
$$K_i=K_i(V)=\{\varphi\in\mbox{End}_F(V)_i:\, \varphi\rho_x=(-1)^{i|x|}\rho_x \varphi \,\,\forall~x\in \mathcal{M}_0\cup \mathcal{M}_1\},\, i=0,1,$$
$$K=K(V)=K_0\oplus K_1,\,\,\mbox{\textit{the supercentralizer} of}\,\,V.$$

The following proposition provides some basic properties of the subsuperalgebras $Z$ and $K$ when we consider irreducible almost faithful modules. 

\begin{pro}[$\cite{ES}$, Proposition $4$]\label{p1}
Assume that $V$ is an irreducible almost faithful module for $\mathcal{M}$. Then
\begin{description}
\item[(i)] $Z_1=0$ and $Z=Z_0$ is an integral domain which acts without torsion on $\mathcal{M}$.
\item[(ii)] $K_0$ is a skew field and any nonzero homogeneous element in $K$ acts bijectively on $V$.
\item[(iii)] The restriction homomorphism $\phi:Z\longrightarrow K_0$($\alpha\longmapsto\alpha|_V$) is one-to-one.
\end{description}
\end{pro}

\subsection{(Super)algebras and modules in $\mathcal{H}$}

In this subsection we define the subvariety $\mathcal{H}$ of the Malcev algebras which was introduced by Filippov and considered in several papers (see for example, \cite{F, F1}) and in the case of superalgebras by Elduque and Shestakov \cite{E, ES, Sh}.

Let $\mathcal{M}$ be a Malcev algebra and consider the subspace $H(\mathcal{M})$ generated by the elements $h(y,z,t,u,x)$; $H(\mathcal{M})$ is an ideal of $\mathcal{M}$ (see $\cite{F}$). 
The variety $\mathcal{H}$ (over $F$) is defined as the class of Malcev algebras $\mathcal{M}$ over $F$ that satisfy the identity $h(y,z,t,u,x)=0$, that is, $H(\mathcal{M})=0$. 

Some elements of the variety $\mathcal{H}$ are the 3-dimensional simple Lie algebra $\mathfrak{s l}(2,F)$ and the 7-dimensional exceptional Malcev algebra $\mathbb{M}$ over its centroid $\Gamma$, which is a field.

If we consider in $\mathcal{M}$ the function $\cite{F1}$
$$p(x,y,z,t)=-\{zt,x,y\}-\{yt,z,x\}+\{xt,y,z\}.$$
We have the following result.
\begin{lem}[$\cite{Sh,ES}$, Lemma $2$]\label{l1}
Let $\mathcal{M}$ be a Malcev algebra and assume that $H(\mathcal{M})=0$. Then for any $x,y,z,t,u\in \mathcal{M}$:
\begin{equation}
    p(x,y,z,t)u=p(xu,y,z,t),
\end{equation}
\begin{equation}
    p(x,y,z,ut)=p(x,u,t,yz).
\end{equation}
\end{lem}

If $\mathcal{M}$ is a Malcev superalgebra. Similarly, we define $\Tilde{p}(x,y,z,t)$. Consider the subspace $\tilde{H}(\mathcal{M})$ generated by the elements $\Tilde{p}(x,y,z,t)$. Without generating confusion we denote again by $\mathcal{H}$ the variety of the class of Malcev superalgebras in which $\Tilde{h}=0$, that is, $\tilde{H}(\mathcal{M})=0$.

Let $\mathcal{M}$ be a Malcev superalgebra in $\mathcal{H}$ and let $V=V_0\oplus V_1$ be a module for $\mathcal{M}$. The module $V$ is called a module for $\mathcal{M}$ in the variety $\mathcal{H}$ or $\mathcal{H}$-\textit{module} if the split extension $E=\mathcal{M}\oplus V$ is again a Malcev superalgebra in $\mathcal{H}$ with a multplicación $(\ref{super1})$.

Elduque studied the subvariety $\mathcal{H}$ of Malcev superalgebras and their representations. He obtained the following results:

\begin{theo}[$\cite{E}$, Theorem $5$]\label{t1}
Let $\mathcal{M}$ be a simple superalgebra in the variety $\mathcal{H}$ and let $V$ be an $\mathcal{H}$-module for $\mathcal{M}$. Then $V$ is completely reducible.
\end{theo}

As a consequence of the proof of the theorem $\ref{t1}$, the author deduced:
\begin{coro}\label{c2}
Let $\mathcal{M}$ be a simple algebra in $\mathcal{H}$, $V$ an $\mathcal{H}$-module for $\mathcal{M}$, and $\Gamma$ the centroid of $\mathcal{M}$. Then $V=N_{V} \oplus J_{V}$ where the submodules $N_{V}$ and $J_{V}$ are given by $N_{V}=\{v \in V: J(v, \mathcal{M}, \mathcal{M})=0\}$ and $J_{V}=J(V, \mathcal{M}, \mathcal{M})$. Moreover,

\begin{description}
\item[(i)] If $\operatorname{dim}_{\Gamma} \mathcal{M}=3$, then $N_{V}=Z_{V} \oplus \bar{N}_{V}$, with $Z_{V}=\{v \in V: v \mathcal{M}=0\}$, $\bar{N}_{V}=V \mathcal{M} \cap N_{V} .$ Besides, $\bar{N}_{V}$ is a direct sum of copies of the adjoint module for $\mathcal{M}$ and $J_{V}$ is a direct sum of copies of the unique irreducible non-Lie module for $\mathcal{M}$.
\item[(ii)] If $\operatorname{dim}_{\Gamma} \mathcal{M}=7$, then $N_{V} \mathcal{M}=0$ and $J_{V}$ is a direct sum of copies of the adjoint module for $\mathcal{M}$.
\end{description}
\end{coro}

Theorem $\ref{t1}$ can be easily extended to the case in which $\mathcal{M}$ is a finite direct sum of simple algebras in $\mathcal{H}$. 

\begin{pro}\label{p2}
Let $\mathcal{M}$ be a finite direct sum of simple algebras in $\mathcal{H}$, $\mathcal{M}=\mathcal{M}_{1} \oplus \cdots \oplus \mathcal{M}_{n}$, with $\mathcal{M}_{i}$ simple in $\mathcal{H}$, and let $V$ be an $\mathcal{H}$-module for $\mathcal{M}.$ Then $V=V_{0} \oplus V_{1} \oplus \cdots \oplus V_{n}$, where $V_{i}=V \mathcal{M}_{i}=V_{i} \mathcal{M}_{i}$ for $i=1, \ldots, n$
and $V_{i} \mathcal{M}_{j}=0$ if $i \neq j$.
\end{pro}

From simple non-Lie Malcev algebras we can suppress the hypothesis of the module being considered in the variety $\mathcal{H}$.

\begin{theo}[$\cite{E}$, Theorem $8$]\label{t2}
Let $\mathcal{M}$ be a finite direct sum of simple non-Lie Malcev algebras. Then any module for $\mathcal{M}$ (in the variety of Malcev algebras) is a $\mathcal{H}$-module.
\end{theo}

Combining Theorems $\ref{t1}$ and $\ref{t2}$, the following result establishes the complete reducibility of Malcev modules over direct sum of simple non-Lie Malcev algebras.

\begin{coro}[$\cite{E}$, Corollary $9$]\label{c1}
Let $\mathcal{M}$ be a finite direct sum of simple non-Lie Malcev algebras. Then any representation of $\mathcal{M}$ is completely reducible.
\end{coro}

As for any variety algebras, the study of the irreducible representations of Malcev (super)algebras is a key point in its investigation. Thus in this direction, Elduque and Shestakov \cite{ES} described the irreducible Non-Lie modules for Malcev superalgebras. 

To finish this subsection we will mention some of the results obtained in \cite{ES} that will be necessary to obtain our main results.

Consider an irreducible $\mathcal{H}$-module $V$ for $\mathcal{M}$ (in the variety of Malcev superalgebras), then $E=\mathcal{M}\oplus V$ is a Malcev superalgebra in $\mathcal{H}$ that is $\tilde{H}(E)=0$. So, by Lemma $\ref{l1}$ we have 
\begin{equation}\label{e5}
    \tilde{p}(x,y,z,t)u=(-1)^{|u|(|y|+|z|+|t|)}\tilde{p}(xu,y,z,t),
\end{equation}
\begin{equation}\label{e6}
    \tilde{p}(x,y,z,ut)=(-1)^{(|u|+|t|)(|y|+|z|)}\tilde{p}(x,u,t,yz)
\end{equation}
for any homogeneous $x,y,z,t,u\in E$. 

Let us define the operator $\tilde{\alpha}(y,z,t)\in\mbox{End}_F(E)$ by
$$x\tilde{\alpha}(y,z,t)=\tilde{p}(x,y,z,t).$$
Then $\tilde{\alpha}$ is super-skewsymmetric on its arguments and it follows from $(\ref{e5})$ that $\tilde{\alpha}(y,z,t)$ belongs to $\widetilde{\Gamma}(E)$.

In the following results, the operator $\tilde{\alpha}(y,z,t)$ plays an important role.
So the theorem shows that the only irreducible almost faithful non-Lie module is the adjoint module for $\mathbb{M}$.

\begin{theo}[$\cite{ES}$, Theorem $2$]\label{t3}
Let $V$ be an irreducible almost faithful non-Lie module for the Malcev superalgebra $\mathcal{M}$ and let $E=\mathcal{M}\oplus V$. Assume $\tilde{H}(E)=0\neq \tilde{\alpha}(V,\mathcal{M},\mathcal{M})$. Then, $\mathcal{M}_{1}=0$, $\mathcal{M}=\mathcal{M}_0$ is a simple non-Lie Malcev algebra, so seven-dimensional over its centroid, and $V$ is the regular (adjoint) module.
\end{theo} 

To describe the structure of every irreducible non-Lie Malcev module $V$ over $\mathbb{M}$, we will do the following: let $v\in V$ and $x,y\in \mathbb{M}$ be with $\tilde{\alpha}(v,x,y)\neq 0$, this induces the isomorphism of modules that was proved in the Theorem $\ref{t3}$
$$\alpha:\mathbb{M}\longrightarrow V$$
given by $\alpha(z)=z\tilde{\alpha}(v,x,y)$. Then $V\cong\mathbb{M}\alpha$ as modules. 


\begin{pro}[$\cite{ES}$, Proposition $8$]\label{p3}
Let $V$ be an irreducible almost faithful non-Lie module for the Malcev superalgebra $\mathcal{M}$ and let $E=V\oplus \mathcal{M}$. If $\tilde{H}(E)=0=\tilde{\alpha}(V,\mathcal{M},\mathcal{M})$. Then, $\mathcal{M}_1=0$, $\mathcal{M}=\mathcal{M}_0$ is a Lie algebra and the map 
$\varphi:\mathcal{M}\longrightarrow\textup{End}_F (V)$, $x\longmapsto -\frac{1}{2}\rho_x$, is a Lie representation of $\mathcal{M}$.
\end{pro}


\section{Factorization Theorems}

In this section we consider Malcev algebras $\mathcal{M}$ in $\mathcal{H}$ containing the $7$-dimensional exceptional Malcev algebra $\mathbb{M}$ such that $m\mathbb{M}\neq 0$ for any $m\neq 0$ from $\mathcal{M}$.
Also, we consider arbitrary Malcev superalgebra $\mathcal{M}=\mathcal{M}_0\oplus \mathcal{M}_1$ in the variety $\mathcal{H}$ whose even part $\mathcal{M}_0$ contains $\mathbb{M}$ with $m\mathbb{M}\neq 0$ for any homogeneous element $0\neq m\in \mathcal{M}_0\cup \mathcal{M}_1$. In the general context, we can drop this assumption and so we describe the Malcev (super)algebras containing $\mathbb{M}$.

\subsection{Factorization Theorem for Malcev algebras}

The objective of this subsection is to prove an analogue of the Kronecker Factorization Theorem for Malcev algebras in $\mathcal{H}$ that contain $\mathbb{M}$.

\begin{theo}\label{t4}
Let $\mathcal{M}$ be a Malcev algebra in $\mathcal{H}$ such that $\mathcal{M}$ contains $\mathbb{M}$, with $m\mathbb{M}\neq 0$ for any $m\neq 0$ from $\mathcal{M}$. Then $\mathcal{M}\cong\mathbb{M}\otimes_F \mathcal{U}$
for a certain commutative associative algebra $\mathcal{U}$.
\end{theo}
\begin{dem}
Consider $\mathcal{M}$ as a Malcev $\mathbb{M}$-module. 
Then by Corollary $\ref{c1}$, $\mathcal{M}$ is completely reducible, that is, $\mathcal{M}=\sum_{i}\oplus V_{i}$, where $V_{i}$ is an irreducible almost faithful non-Lie Malcev module over $\mathbb{M}$, that is, $V_i\cong \textup{Reg}_{i}\,\mathbb{M}\cong \textup{Reg}\,\mathbb{M}$ for all $i$. 
Also by Proposition \ref{p3} we have $0\neq \alpha_i (V_i,\mathbb{M},\mathbb{M})$ because $\mathbb{M}$ is a non-Lie Malcev algebra. Let us take elements $v_{i}\in V_i$ and $a_{i},b_{i}\in\mathbb{M}$ such that
$$\alpha_{i}=\alpha_{i}(v_{i},a_{i},b_{i})\neq 0.$$
Thus 
\begin{equation*}
    \mathcal{M}=\sum_{i}\oplus \mathbb{M}\alpha_{i}
\end{equation*} 
where for each $i$ we have $\textup{Reg}_{i}\,\mathbb{M}\cong\mathbb{M}\alpha_{i}$.
Moreover, since $\mathcal{M}$ satisfies $H(\mathcal{M})=0,$ then
by Lemma $\ref{l1}$ for all $x,y\in \mathcal{M}$ we have
\begin{equation*}
    p(x,v_{i},a_{i},b_{i})y=p(xy,v_{i},a_{i},b_{i}),
\end{equation*}
$(x\alpha_{i})y=(xy)\alpha_{i}$; thus $\alpha_{i}\in\Gamma(\mathcal{M})$, where $\Gamma(\mathcal{M})$ is \textit{the centroid} of $\mathcal{M}$.

Let $\mathcal{U}=\sum_{i} F\alpha_i$ denote the span of all $\alpha_i$. Then $\mathcal{M}=\mathbb{M}\,\mathcal{U}.$ So let's go to prove that $\mathcal{U}$ is a subalgebra of $\Gamma(\mathcal{M})$. 

First we consider $\mathbb{M}=\mathcal{M}_7(F)$, so we fix arbitrary elements $\alpha, \beta\in\mathcal{U}$ and $a,b\in\mathcal{M}_7(F)$. So $a\alpha\beta\in \mathcal{M}$, then $a\alpha\beta=\sum_i a_i\alpha_i$ for some $a_i\in\mathcal{M}_7(F)$ and $\alpha_i\in \mathcal{U}$. We denote $w=a\alpha\beta-\sum_i a_i\alpha_i=0$
\begin{equation}\label{e7}
    w=(\sum^{7}_{t=1}\lambda_t e_t)\alpha\beta-\sum_{i}(\sum^{7}_{t=1}\lambda_{ti} e_t)\alpha_i=0,
\end{equation}
where $a=\sum^{7}_{t=1}\lambda_t e_t$, $a_i=\sum^{7}_{t=1}\lambda_{ti} e_t$ and $0\neq\lambda_t,\lambda_{ti}\in F$, $t=1,2,\dots,7$.
As $\alpha\beta\in \Gamma(\mathcal{M})$ and using the multiplication table of $\mathcal{M}_7(F)$, we have 
$$0=e_2 w=(-2\lambda_1 e_2+2\lambda_3 e_7-2\lambda_4 e_6+\lambda_5 e_1)\alpha\beta-
\sum_{i}(-2\lambda_{1i} e_2+2\lambda_{3i} e_7-2\lambda_{4i} e_6+\lambda_{5i} e_1)\alpha_i,$$
\begin{equation*}
    0=e_2(e_2 w)=-2(\lambda_5 e_2)\alpha\beta+2\sum_{i}(\lambda_{5i} e_2)\alpha_i,
\end{equation*}
which implies, $e_2\alpha\beta=e_2\widetilde{\alpha}$, where $\widetilde{\alpha}=\sum_{i}(\lambda^{-1}_5\lambda_{5i})\alpha_i\in\mathcal{U}$. Hence, it is easy to see $e_t(\alpha\beta-\widetilde{\alpha})=0$ for $t=1,2,\dots,7$; so
$$\mathcal{M}_7(F)(\alpha\beta-\widetilde{\alpha})=0$$
and $\alpha\beta-\widetilde{\alpha}=0$. Therefore $\alpha\beta=\widetilde{\alpha}\in\mathcal{U}$; $\mathcal{U}\mathcal{U}\subseteq\mathcal{U}$.

Now let $m\in \mathcal{M}$ and $a\in \mathcal{M}_7(F)$ be, then
\begin{equation*}
\begin{split}
((m\alpha)\beta)a&=(m\alpha)(a\beta)=((m\alpha)a)\beta=(m(a\alpha))\beta\\
                    &=(m\beta)(a\alpha)=((m\beta)a)\alpha=((m\beta)\alpha)a.
\end{split}
\end{equation*}
If $[\alpha,\beta]=\alpha\beta-\beta\alpha$, we have $$(\mathcal{M}[\alpha,\beta])\mathcal{M}_7(F)=0,$$ so $\mathcal{M}[\alpha,\beta]=0$ and $[\alpha,\beta]|_{\mathcal{M}}=0$. In particular, $[\alpha,\beta]|_{V_{i}}=0$ for any irreducible component $V_i$ of $\mathcal{M}$, then by Proposition $\ref{p1}$(iii) $[\alpha,\beta]=0$ because $\phi:Z\longrightarrow K$($\alpha\longmapsto\alpha|_{V_{i}}$) is one-to-one. Therefore $[\mathcal{U},\mathcal{U}]=0$; hence $\mathcal{U}$ is a commutative and associative algebra.

Also
\begin{equation*}
\begin{split}
(a\alpha)(b\beta)&=(a(b\beta))\alpha=((ab)\beta)\alpha\\
                    &=(ab)(\beta\alpha)=(ab)(\alpha\beta).
\end{split}
\end{equation*}

Let $v=\sum^{7}_{i=1}e_i \alpha_i=0$ be, where $\alpha_i\in \mathcal{U}$. Then
$0=e_4 v=-2e_4 \alpha_1+2e_6 \alpha_2-2e_5 \alpha_3$ and
$0=e_5(e_4 v)=-4e_4 \alpha_2.$
Hence $e_4 \alpha_2=0$. Also
\begin{equation*}
\begin{split}
0=e_2(e_4 \alpha_2)&=(e_2e_4)\alpha_2=-2e_6 \alpha_2,\\
0=e_3(e_4 \alpha_2)&=(e_3e_4)\alpha_2=2e_5 \alpha_2,\\
0=e_7(e_4 \alpha_2)&=(e_7e_4)\alpha_2=-e_1 \alpha_2,
\end{split}
\end{equation*}
so $e_6 \alpha_2=e_5 \alpha_2=e_1 \alpha_2=0$, and
\begin{equation*}
\begin{split}
0=e_2(e_1 \alpha_2)&=(e_2e_1)\alpha_2=-2e_2 \alpha_2,\\
0=e_3(e_1 \alpha_2)&=(e_3e_1)\alpha_2=-2e_3 \alpha_2,\\
0=e_7(e_1 \alpha_2)&=(e_7e_1)\alpha_2=2e_7 \alpha_2,
\end{split}
\end{equation*}
so $e_2 \alpha_2=e_3 \alpha_2=e_7 \alpha_2=0$. Thus $\mathcal{M}_7(F)\alpha_2=0$ and $\alpha_2=0$. Similarly $\alpha_i=0$ for $i=1,3,\dots,7$. Therefore $\mathcal{M}\cong\mathcal{M}_7(F)\otimes_F \mathcal{U}$.

If $\mathbb{M}=\mathbb{M}_7$. Without generating confusion we denote again by $\mathcal{U}$ the vector space generated by all the isomorphisms $\alpha_{i}$ between the irreducible $\mathbb{M}_7$-modules and $\mbox{Reg}\,\mathbb{M}_7$.

As in above, we will prove that $\mathcal{U}$ is a subalgebra of $\Gamma(\mathcal{M})$. We fix arbitrary elements $\alpha, \beta\in\mathcal{U}$ and $a,b\in\mathbb{M}_7$. So $a\alpha\beta\in \mathcal{M}$, then $a\alpha\beta=\sum_t a_t\alpha_t$ for some $a_t\in\mathbb{M}_7$ and $\alpha_t\in \mathcal{U}$. Denote $w=a\alpha\beta-\sum_t a_t\alpha_t=0$, that is
\begin{equation*}
\begin{split}
w&=(\lambda_{1}i+\lambda_{2}j+\lambda_{3}k+\lambda_{4}v+\lambda_{5}vi+\lambda_{6}vj+\lambda_{7}vk)\alpha\beta\\
&\,\,\,\,\,\,\,-\sum_{t}(\lambda_{1t}i+\lambda_{2t}j+\lambda_{3t}k+\lambda_{4t}v+\lambda_{5t}vi+\lambda_{6t}vj+\lambda_{7t}vk)\alpha_t=0,
\end{split}
\end{equation*}
where $a=\lambda_{1}i+\lambda_{2}j+\lambda_{3}k+\lambda_{4}v+\lambda_{5}vi+\lambda_{6}vj+\lambda_{7}vk$, $a_t=\lambda_{1t}i+\lambda_{2t}j+\lambda_{3t}k+\lambda_{4t}v+\lambda_{5t}vi+\lambda_{6t}vj+\lambda_{7t}vk$ and $0\neq\lambda_{s}, \lambda_{st}\in F$, $s=1,\dots, 7$. 
So, using the multiplication table of $\mathbb{M}_7$
\begin{equation*}
\begin{split}
0=iw&=(2\lambda_{2}k-2\lambda_{3}j-2\lambda_{4}vi+2\lambda_{5}v-2\lambda_{6}vk+2\lambda_{7}vj)\alpha\beta\\
&\,\,\,\,\,\,\,-\sum_{t}(2\lambda_{2t}k-2\lambda_{3t}j-2\lambda_{4t}vi+2\lambda_{5t}v-2\lambda_{6t}vk+2\lambda_{7t}vj)\alpha_t,
\end{split}
\end{equation*}
\begin{equation*}
\begin{split}
0=k(iw)&=(4\lambda_{3}i+4\lambda_{4}vj-4\lambda_{5}vk-4\lambda_{6}v+4\lambda_{7}vi)\alpha\beta\\
&\,\,\,\,\,\,\,-\sum_{t}(4\lambda_{3t}i+4\lambda_{4t}vj-4\lambda_{5t}vk-4\lambda_{6t}v+4\lambda_{7t}vi)\alpha_t,
\end{split}
\end{equation*}
\begin{equation*}
\begin{split}
0=i(k(iw))&=(-8\lambda_{4}vk-8\lambda_{5}vj+8\lambda_{6}vi+8\lambda_{7}v)\alpha\beta\\
&\,\,\,\,\,\,\,-\sum_{t}(-8\lambda_{4t}vk-8\lambda_{5t}vj+8\lambda_{6t}vi+8\lambda_{7t}v)\alpha_t,
\end{split}
\end{equation*}
\begin{equation*}
\begin{split}
0=(vk)(i(k(iw)))&=(16\lambda_{5}\gamma i+16\lambda_{6}\gamma j-16\lambda_{7}\gamma k)\alpha\beta\\
&\,\,\,\,\,\,\,-\sum_{t}(16\lambda_{5t}\gamma i+16\lambda_{6t}\gamma j-16\lambda_{7t}\gamma k)\alpha_t,
\end{split}
\end{equation*}
\begin{equation*}
0=i((vk)(i(k(iw))))=(32\lambda_{6}\gamma k+32\lambda_{7}\gamma j)\alpha\beta -\sum_{t}(32\lambda_{6t}\gamma k+32\lambda_{7t}\gamma j)\alpha_t,
\end{equation*}
\begin{equation*}
0=k(i((vk)(i(k(iw)))))=(-64\lambda_{7}\gamma i)\alpha\beta -\sum_{t}(-64\lambda_{7t}\gamma i)\alpha_t.
\end{equation*}
As $\lambda_{7}, \gamma\neq 0$, we have
\begin{center}
    $i(\alpha\beta-\hat{\alpha})=0$
\end{center}
where $\hat{\alpha}=\sum_{t}(\lambda^{-1}_7\lambda_{7t})\alpha_t\in\mathcal{U}$. Hence, it is easy to see $x(\alpha\beta-\hat{\alpha})=0$ for any base elements $x\in\left\{i,j,k,v,vi,vj,vk\right\}$; hence 
$$\mathbb{M}_7(\alpha\beta-\hat{\alpha})=0$$
and $\alpha\beta-\hat{\alpha}=0$. Thus $\alpha\beta=\hat{\alpha}\in\mathcal{U}$; $\mathcal{U}\mathcal{U}\subseteq\mathcal{U}$.

Also, following the same ideas of the case $\mathcal{M}_7(F)$, we can prove that $[\mathcal{U},\mathcal{U}]=0$; thus the algebra $\mathcal{U}$ is associative and commutative. Also
$
(a\alpha)(b\beta)=(ab)(\alpha\beta)
$
for any $\alpha,\beta\in \mathcal{U}$ and $a,b\in\mathbb{M}_7$.

Finally, we shall prove that $\mathcal{U}$ is free over $\mathbb{M}_7$. We consider 
$\tilde{w}=i\alpha_1+j\alpha_2+k\alpha_3+v\alpha_4 +(vi)\alpha_5+(vj)\alpha_6+(vk)\alpha_7=0$, where $\alpha_i\in \mathcal{U}$, $i=1,2,\dots,7$. Then, by the multiplication table of $\mathbb{M}_7$
\begin{center}
$0=i\tilde{w}=2k\alpha_2-2j\alpha_3-2(vi)\alpha_4+2v\alpha_5-2(vk)\alpha_6+2(vj)\alpha_7,$
\end{center}
\begin{center}
$0=k(i\tilde{w})=4i\alpha_3+4(vj)\alpha_4-4(vk)\alpha_5-4v\alpha_6+4(vi)\alpha_7,$
\end{center}
\begin{center}
$0=i(k(i\tilde{w}))=-8(vk)\alpha_4-8(vj)\alpha_5+8(vi)\alpha_6+8v\alpha_7,$
\end{center}
\begin{center}
$0=(vk)(i(k(i\tilde{w})))=16\gamma i\alpha_5+16\gamma j\alpha_6-16\gamma k\alpha_7,$
\end{center}
\begin{center}
$0=i((vk)(i(k(i\tilde{w}))))=32\gamma k\alpha_6+32\gamma j\alpha_7,$
\end{center}
\begin{center}
$0=k(i((vk)(i(k(i\tilde{w})))))=-64\gamma i\alpha_7.$
\end{center}
As $\gamma\neq 0$, of the last equality we have 
$0=i\alpha_7.$ Hence, it is easy to see $x\alpha_7=0$ for any base elements $x\in\left\{i,j,k,v,vi,vj,vk\right\}$; hence 
$$\mathbb{M}_7\alpha_7=0$$
and $\alpha_7=0$. Similarly $\alpha_i=0$ for $i=1,2,3,\dots,6$. Therefore $\mathcal{M}\cong\mathbb{M}_7\otimes_{F} \mathcal{U}$.

The theorem is proved.
\end{dem}

Using  Theorem $\ref{t4}$ we  can  easily  drop  the  assumption  that $m\mathbb{M}\neq 0$ for any $m\neq 0$ from $\mathcal{M}$.

\begin{theo}\label{t7}
Let $\mathcal{M}$ be a Malcev algebra in 
$\mathcal{H}$ containing $\mathbb{M}$. Then $\mathcal{M}$ is isomorphic to a direct sum $(\mathbb{M}\otimes_{F} A)\oplus \mathcal{N}$, where $A$ is a certain commutative associative algebra.
\end{theo}
\begin{dem}
Consider $\mathcal{M}$ as a Malcev $\mathbb{M}$-module. Then by the proof of Theorem $\ref{t4}$ we have the first decomposition $\mathcal{M}=\sum_{i}\oplus \mathbb{M}\alpha_{i}$, where $\alpha_{i}$ belong to $\Gamma(\mathcal{M})$ (the centroid of $\mathcal{M}$). From this decomposition it follows that $\mathcal{M}=\mathbb{M}\,\mathcal{U}$, where $\mathcal{U}=\sum_{i} F\alpha_i$ denote the span of all $\alpha_i$. 

On the other hand, from Corollary $\ref{c2}$,
$\mathcal{M}=N_{\mathcal{M}}\oplus J_{\mathcal{M}}$, where the submodules $N_{\mathcal{M}}$ and $J_{\mathcal{M}}$ are given by $N_{\mathcal{M}}=\{x \in \mathcal{M}: J(x, \mathbb{M}, \mathbb{M})=0\}$ and $J_{\mathcal{M}}=J(\mathcal{M}, \mathbb{M}, \mathbb{M})$. Moreover, by $\textup{(ii)}$
of Corollary $\ref{c2}$,
$N_{\mathcal{M}}\mathbb{M}=0$
and $J_{\mathcal{M}}$ is a direct sum of the copies of the adjoint module of $\mathbb{M}$, that is, $J_{\mathcal{M}}=\sum_{j}\oplus \textup{Reg}_{j}\,\mathbb{M}$, where $\textup{Reg}_{j}\,\mathbb{M}\cong \textup{Reg}\,\mathbb{M}$ for all $j$. So 
$J_{\mathcal{M}}=\sum_{j}\oplus \mathbb{M}\alpha_{j}$, where $\alpha_{j}\in\Gamma(\mathcal{M})$ and they are some of the elements that appear in the first decomposition of $\mathcal{M}$. Denote $A=\sum_{j} F\alpha_j$ the span of all $\alpha_j$. Then $J_{\mathcal{M}}=\mathbb{M} A.$

For arbitrary elements $x\in N_{\mathcal{M}}$, $a,b,c\in\mathbb{M}$, $\alpha\in\mathcal{U}$, the equality $\mathcal{M}=\mathbb{M}\,\mathcal{U}$ and using $(\ref{e8})$ we have
\begin{center}
    $J(x(a\alpha),b,c)=J(xa,b,c)\alpha=(xJ(a,b,c)+J(x,b,c)a-2J(x,a,bc))\alpha=0;$
\end{center}
so $N_{\mathcal{M}}$ is an ideal of $\mathcal{M}$. Now as $A^{2}\subseteq\Gamma(\mathcal{M})$ and $\mathbb{M}\,A^{2}\subseteq N_{\mathcal{M}} + J_{\mathcal{M}}$ we obtain
\begin{center}
    $J^{2}_{\mathcal{M}}=(\mathbb{M}\,A)(\mathbb{M}\,A)=\mathbb{M}^{2}\,A^{2}=(\mathbb{M}\,A^{2})\mathbb{M}\subseteq N_{\mathcal{M}}\mathbb{M} + J_{\mathcal{M}}\mathbb{M}\subseteq J_{\mathcal{M}}$,
\end{center}
and $J^{2}_{\mathcal{M}}\subseteq J_{\mathcal{M}}$; so $J_{\mathcal{M}}$ is a subalgebra of $\mathcal{M}$. Then by Theorem $\ref{t4}$, $J_{\mathcal{M}}\cong\mathbb{M}\otimes_{F} A$, where $A$ is a commutative associative algebra. Denote $\mathcal{N}:=N_{\mathcal{M}}\triangleleft\mathcal{M}$. Thus 
\begin{center}
    $\mathcal{M}\cong(\mathbb{M}\otimes_{F} A)\oplus \mathcal{N}.$
\end{center}

The theorem is proved.
\end{dem}


\subsection{Factorization Theorem for Malcev superalgebras}

The objective of this subsection is to prove an analogue of the Kronecker Factorization Theorem for Malcev superalgebras in $\mathcal{H}$ whose even part contains the 7-dimensional exceptional Malcev algebra $\mathbb{M}$. 

\begin{theo}\label{t5}
Let $\mathcal{M}=\mathcal{M}_0\oplus \mathcal{M}_1$ be a Malcev superalgebra in $\mathcal{H}$ such that $\mathcal{M}_0$ contains $\mathbb{M}$, with $m\mathbb{M}\neq 0$ for any homogeneous element $0\neq m\in \mathcal{M}_0\cup \mathcal{M}_1$. Then $\mathcal{M}\cong\mathbb{M}\otimes_{F} U$ for a certain supercommutative associative superalgebra $U$.
\end{theo}
\begin{dem}
Consider $\mathcal{M}_0$ as a Malcev $\mathbb{M}$-module. Then, as in the Theorem $\ref{t4}$, every irreducible non-Lie Malcev $\mathbb{M}$-module is almost faithful and so the regular module $\mbox{Reg}\,\mathbb{M}$.

As $\mathcal{M}_0$ is a Malcev module over $\mathbb{M}$, by Corollary $\ref{c1}$, $\mathcal{M}_0$ is completely reducible, that is, $\mathcal{M}_0=\sum_i\oplus V_i$, where $V_i$ is an irreducible almost faithful non-Lie Malcev $\mathbb{M}$-module. 
Again by Proposition \ref{p3} we have $0\neq \alpha_i (V_i,\mathbb{M},\mathbb{M})$ because $\mathbb{M}$ is a non-Lie Malcev algebra. Let us take elements $v_{i}\in V_i$ and $a_{i},b_{i}\in\mathbb{M}$ such that 
$$\alpha_{i}=\alpha_{i}(v_{i},a_{i},b_{i})\neq 0$$
since $\mathbb{M}$ is a non-Lie Malcev algebra.
Therefore
\begin{equation*}
    \mathcal{M}_0=\sum_{i}\oplus \mathbb{M}\alpha_i
\end{equation*}
because for each $i$ we have $V_i\cong\mathbb{M}\alpha_i$. As in the Theorem $\ref{t4}$, $\alpha_{i}\in\widetilde{\Gamma}_0 (\mathcal{M})$.

Let $U_0=\sum_{i} F\alpha_i$ denote the span of all $\alpha_i$. Then $\mathcal{M}_0=\mathbb{M}U_0$.

Also we can consider $M_1$ as a Malcev module over $\mathbb{M}$, by Corollary $\ref{c1}$, $\mathcal{M}_1$ is completely reducible, that is, $\mathcal{M}_1=\sum_j\oplus W_j$, where $W_j$ is an irreducible almost faithful non-Lie Malcev $\mathbb{M}$-module. Thus
\begin{equation*}
    \mathcal{M}_1=\sum_j\oplus \mathbb{M}\beta_j
\end{equation*}
because for each $j$ we have $W_j\cong\mathbb{M}\beta_j$, with
\begin{equation*}
    \beta_j=\beta_j(w_j,a_j,b_j)\neq 0
\end{equation*}
for some $w_j\in W_j$ and $a_j,b_j\in\mathbb{M}$ because by Proposition \ref{p3} we have $0\neq \beta_j (W_j,\mathbb{M},\mathbb{M})$ since $\mathbb{M}$ is a non-Lie Malcev algebra. Also $\beta_{j}\in\widetilde{\Gamma}_1 (\mathcal{M})$. 

We denote by $U_1=\sum_{j} F\beta_j$ the span of all $\beta_j$. Then $\mathcal{M}_1=\mathbb{M}U_1$.

We have 
\begin{equation*}
    \mathcal{M}=\mathbb{M} \,U_0\oplus\mathbb{M}\,U_1.
\end{equation*}
Let's go to prove that $U=U_0\oplus U_1$ is a subsuperalgebra of $\widetilde{\Gamma}(\mathcal{M})$. For all $a,b\in\mathbb{M}$ we have
\begin{equation}\label{cen1}
(a\alpha)(b\beta)=(-1)^{|\alpha|(|b|+|\beta|)}(a(b\beta))\alpha=(-1)^{|\alpha|(|b|+|\beta|)}((ab)\beta)\alpha
                    =(-1)^{|\alpha|(|b|+|\beta|)}(ab)(\beta\alpha).
\end{equation}
So by $(\ref{cen1})$ and using the relationships 
\begin{equation*}
    M_iM_j\subseteq M_{(i+j)mod\,2},\,\,\,\,i,j=0,1
\end{equation*}
of the $\mathbb{Z}_2$-grading of $\mathcal{M}$ and the proof of the Theorem $\ref{t4}$, we get
$$U_i U_j\subseteq U_{(i+j)mod\,2},\,\,\,\,i,j=0,1.$$
Hence $U=U_0\oplus U_1$ is a subsuperalgebra.

We fix arbitrary homogeneous elements $\alpha, \beta\in U$, $m\in \mathcal{M}$ and $a\in\mathbb{M}$, then
\begin{equation*}
\begin{split}
((m\alpha)\beta)a&=(-1)^{|a||\beta|}(m\alpha)(a\beta)\\
                 &=(-1)^{|a|(|\beta|+|\alpha|)+|\alpha||\beta|}(m(a\alpha))\beta\\
                    &=(-1)^{|\alpha|(|a|+|\beta|)}((m\beta)a)\alpha\\
                    &=(-1)^{|\alpha||\beta|}((m\beta)\alpha)a.
\end{split}
\end{equation*}
If $[\alpha,\beta]=\alpha\beta-(-1)^{|\alpha||\beta|}\beta\alpha$, we have 
$$(\mathcal{M}[\alpha,\beta])\mathbb{M}=0,$$
so $[\alpha,\beta]|_{\mathcal{M}}=0$. In particular, $[\alpha,\beta]|_{V_{i}}=0$ or $[\alpha,\beta]|_{W_{j}}=0$ for any irreducible components $V_i$, $W_j$ of $\mathcal{M}$, then by Proposition $\ref{p1}$(iii) $[\alpha,\beta]=0$ because $\phi:Z\longrightarrow U$($\alpha\longmapsto\alpha|_{V_{i}}$) or $\varphi:Z\longrightarrow U$($\beta\longmapsto\beta|_{W_{j}}$) are one-to-one. Therefore $[U,U]=0$; hence $U$ is a supercommutative and associative superalgebra.

Also, by $(\ref{cen1})$ and $[U,U]=0$
\begin{equation*}
    (a\alpha)(b\beta)=(-1)^{|\alpha|(|b|+|\beta|)}(ab)(\beta\alpha)=(-1)^{|\alpha|(|b|+|\beta|)}(-1)^{|\alpha||\beta|}(ab)(\alpha\beta)
                    =(-1)^{|\alpha||b|}(ab)(\alpha\beta).
\end{equation*}
for all $a,b\in\mathbb{M}$.

As in the Theorem $\ref{t4}$, using the multiplication table of $\mathbb{M}_7$ and $\mathcal{M}_7(F)$, it is easy to see that the superalgebra $U$ is free over $\mathbb{M}$.
Therefore $\mathcal{M}\cong\mathbb{M}\otimes_{F} U$.

The theorem is proved.
\end{dem}

Using  Theorem $\ref{t5}$  we  can  drop  the  assumption  that $m\mathbb{M}\neq 0$ for any homogeneous element $0\neq m\in \mathcal{M}_0\cup \mathcal{M}_1$. The proof of the Theorem $\ref{t8}$ mimics the Theorem $\ref{t7}$ proof of the corresponding result for Malcev algebras in $\mathcal{H}$ that contains the $7$-dimensional exceptional Malcev algebra.

\begin{theo}\label{t8}
Let $\mathcal{M}=\mathcal{M}_0\oplus \mathcal{M}_1$ be a Malcev superalgebra in $\mathcal{H}$ such that $\mathcal{M}_0$ contains $\mathbb{M}$. Then $\mathcal{M}$ is isomorphic to a direct sum $(\mathbb{M}\otimes_{F} \mathcal{A})\oplus \mathcal{N'}$, where $\mathcal{A}$ is a certain supercommutative associative superalgebra.
\end{theo}
\begin{dem}
Consider $\mathcal{M}$ as a Malcev $\mathbb{M}$-module. Then by the proof of Theorem $\ref{t5}$ we have the first decomposition $\mathcal{M}=(\sum_{i}\oplus \mathbb{M}\alpha_{i})\oplus (\sum_{j}\oplus \mathbb{M}\beta_{j})$, where $\alpha_{i}\in \widetilde{\Gamma}_{0}(\mathcal{M})$ and $\beta_{j}\in \widetilde{\Gamma}_{1}(\mathcal{M})$. From this decomposition it follows that $\mathcal{M}=\mathbb{M}U$, where $U=U_{0}\oplus U_{1}=\sum_{i} F\alpha_i+\sum_{j} F\beta_j$ is a subsuperalgebra of $\widetilde{\Gamma}(\mathcal{M})$ (the supercentroid of $\mathcal{M}$) spanned by all $\alpha_i$ and $\beta_j$. 

Also from Corollary $\ref{c2}$ we have the decomposition
$\mathcal{M}=N_{\mathcal{M}}\oplus J_{\mathcal{M}}$, where the submodules $N_{\mathcal{M}}$ and $J_{\mathcal{M}}$ are given by $N_{\mathcal{M}}=\{x \in \mathcal{M}: J(x, \mathbb{M}, \mathbb{M})=0\}$ and $J_{\mathcal{M}}=J(\mathcal{M}, \mathbb{M}, \mathbb{M})$. Moreover, by $\textup{(ii)}$
of Corollary $\ref{c2}$
$N_{\mathcal{M}}\mathbb{M}=0$
and $J_{\mathcal{M}}$ is a direct sum of the copies of the adjoint module of $\mathbb{M}$, that is, $J_{\mathcal{M}}=(\sum_{r}\oplus \textup{Reg}_{r}\,\mathbb{M})\oplus (\sum_{s}\oplus \textup{Reg}_{s}\,\mathbb{M})$, where $\textup{Reg}_{r}\,\mathbb{M}\cong \textup{Reg}\,\mathbb{M}$ and $\textup{Reg}_{s}\,\mathbb{M}\cong \textup{Reg}\,\mathbb{M}$ for all $r,s$. So 
$J_{\mathcal{M}}=(\sum_{r}\oplus \mathbb{M}\alpha_{r})\oplus (\sum_{s}\oplus \mathbb{M}\alpha_{s})$, where $\alpha_{r}\in\widetilde{\Gamma}_{0}(\mathcal{M})$ and $\beta_{s}\in\widetilde{\Gamma}_{1}(\mathcal{M})$ and they are some of the elements that appear in the first decomposition of $\mathcal{M}$. Denote $\algA=\sum_{r} F\alpha_r+ \sum_{s} F\alpha_s$ the span of all $\alpha_r$ and $\beta_s$. Then $J_{\mathcal{M}}=\mathbb{M} \algA.$

For arbitrary elements $x\in N_{\mathcal{M}}$, $a,b,c\in\mathbb{M}$, $\alpha\in U_{0}\cup U_{1}$, the equality $\mathcal{M}=\mathbb{M}U$ and the superlinearization of $(\ref{e8})$ imply
\begin{center}
    $J(x(a\alpha),b,c)=\pm J(xa,b,c)\alpha=\pm(\pm xJ(a,b,c)\pm J(x,b,c)a\pm 2J(x,a,bc))\alpha=0;$
\end{center}
so $N_{\mathcal{M}}$ is a superideal of $\mathcal{M}$. Also as  $\algA^{2}\subseteq\widetilde{\Gamma}(\mathcal{M})$ and $\mathbb{M}\,\algA^{2}\subseteq N_{\mathcal{M}} + J_{\mathcal{M}}$ we obtain
\begin{center}
    $J^{2}_{\mathcal{M}}=(\mathbb{M}\algA)(\mathbb{M}\algA)=\mathbb{M}^{2}\algA^{2}=(\mathbb{M}\algA^{2})\mathbb{M}\subseteq N_{\mathcal{M}}\mathbb{M} + J_{\mathcal{M}}\mathbb{M}\subseteq J_{\mathcal{M}}$,
\end{center}
and $J^{2}_{\mathcal{M}}\subseteq J_{\mathcal{M}}$; then $J_{\mathcal{M}}$ is a subsuperalgebra of $\mathcal{M}$ and by Theorem $\ref{t5}$, $J_{\mathcal{M}}\cong\mathbb{M}\otimes_{F} \algA$, where $\algA$ is a supercommutative associative superalgebra. Denote $\mathcal{N}':=N_{\mathcal{M}}\triangleleft\mathcal{M}$. Thus 
\begin{center}
    $\mathcal{M}\cong(\mathbb{M}\otimes_{F} \algA)\oplus \mathcal{N}'.$
\end{center}

The theorem is proved.
\end{dem}

\section{Some equivalence of categories}

In this section, we provide some equivalence of categories. In the case of Malcev algebras, observe that when $\mathbb{M}=\mathcal{M}_7(F)$, the Theorem \ref{t4} states that the Malcev algebra $\mathcal{M}$ in $\mathcal{H}$ is coordinated by $\mathcal{U}$, that is, $\mathcal{M}\cong\mathcal{M}_7 (\mathcal{U})$. Thus the following result establishes that in this case $\mathcal{M}$ is Morita equivalent to $\mathcal{U}$. 

Let $\mathcal{U}$ be a unital associative commutative algebra, and let $C(\mathcal{U})$ be the Cayley-Dickson matrix algebra over $\mathcal{U}$. So introducing the commutator $[A,B]=AB-BA$ in $C(\mathcal{U})$ we obtain the Malcev algebra $C(\mathcal{U})^{(-)}$ and the non-Lie Malcev algebra
$\mathcal{M}_7(\mathcal{U})=C(\mathcal{U})^{(-)}/F\cdot e_{0}$, where $e_{0}=
\begin{bmatrix}
1&0\\
0&1\\
\end{bmatrix}$.
It's clear that $\mathcal{M}_7(\mathcal{U})$ contains the 7-dimensional exceptional Malcev algebra $\mathcal{M}_7(F)$.

Then consider the following categories:
\begin{itemize}
\item let $\mathcal{H}-\textbf{\mbox{Mod}}_\textbf{\mbox{Malc}}\,\mathcal{M}_7(\mathcal{U})$ denote the category of $\mathcal{H}$-modules $V$ over $\mathcal{M}_7(\mathcal{U})$
such that $v\mathcal{M}_7(F)\neq 0$ for any $0\neq v\in V$.
\item let $\textbf{\mbox{Mod}}_\textbf{\mbox{Assoc-Com}}\,\mathcal{U}$ denote the category of unital associative and commutative modules $W$ over $\mathcal{U}$ such that $w\mathcal{U}\neq 0$ for any $0\neq w\in W$.
\end{itemize}
We define a functor $\mathbf{T}$ from the category $\textbf{\mbox{Mod}}_\textbf{\mbox{Assoc-Com}}\,\mathcal{U}$
into the category $\mathcal{H}-\textbf{\mbox{Mod}}_\textbf{\mbox{Malc}}\,\mathcal{M}_7(\mathcal{U}).$

Let $W\in\mbox{Obj}(\textbf{\mbox{Mod}}_\textbf{\mbox{Assoc-Com}}\,\mathcal{U})$ and $E=\mathcal{U}+W$ be the split null extension of $\mathcal{U}$ by $W$ such that $w\mathcal{U}\neq 0$ for any $0\neq w\in W$, $E$ is a unital associative and commutative algebra. Consider the Cayley-Dickson matrix algebra $C(E)$ that contains $C(\mathcal{U})$ as a subalgebra. Also $C(E)$ contains the ideal $C(W)$ which is the set of matrices of $C(E)$ whose entries are in the ideal $W$ of $E$.
So introducing the commutator in $C(E)$ we obtain the non-Lie Malcev algebra 
$$\mathcal{M}_7 (E)=C(E)^{(-)}/F\cdot e_{0}.$$
that contains the non-Lie Malcev algebra $\mathcal{M}_7 (\mathcal{U})=C(\mathcal{U})^{(-)}/F\cdot e_{0}$.
Moreover, since $\mathcal{M}_7(F)$ belongs to $\mathcal{H}$ and $\mathcal{M}_7(E)=\mathcal{M}_7(F)\otimes_{F} E, $ then $\mathcal{M}_7(E)$ belongs to $\mathcal{H}$ and contains the ideal $V=\mathcal{M}_7(W)$.
Therefore, $V=\mathcal{M}_7(W)$ is a $\mathcal{H}$-module over $\mathcal{M}_7(\mathcal{U})$ relative to the multiplication defined in $\mathcal{M}_7(E)$.
Note that $v\mathcal{M}_7(F)\neq 0$ for any $0\neq v\in V$ because $W$ is unital.
We call to $V$ the $\mathcal{M}_7(\mathcal{U})$-module \textbf{associated} with the given module $W$ of $\mathcal{U}$ and denote 
$$V=\mathbf{T}(W).$$
As $E=\mathcal{U}\oplus W$ we have $\mathcal{M}_7(E)=\mathcal{M}_7(\mathcal{U})\oplus V$. Also $W^2=0$ in $E$ which implies that $V^2=0$ in $\mathcal{M}_7(E)$, thus $\mathcal{M}_7(E)$ is the split null extension of $\mathcal{M}_7(\mathcal{U})$ by its module $V$.

We can easily verify that 
$$
\mathbf{T}:\textbf{\mbox{Mod}}_\textbf{\mbox{Assoc-Com}}\,\mathcal{U}\longrightarrow \mathcal{H}-\textbf{\mbox{Mod}}_\textbf{\mbox{Malc}}\,\mathcal{M}_7(\mathcal{U})
$$
is really a functor.
In addition each pair of objects $W$ and $W'$ of $\textbf{\mbox{Mod}}_\textbf{\mbox{Assoc-Com}}\,\mathcal{U}$, the following equality is valid:
$$\mathbf{T}(\mbox{Hom}(W,W'))=\mbox{Hom}(\mathbf{T}(W),\mathbf{T}(W')).$$
Thus $W$ and $W'$ are isomorphic if and only if $\mathbf{T}(W)$ and $\mathbf{T}(W')$
are too.

Similarly the functor $\mathbf{T}$ gives a lattice isomorphism of the lattice of
submodules of $W$ relative to $\mathcal{U}$ onto the lattice of submodules
of $V=\mathbf{T}(W)$ relative to $\mathcal{M}_7(\mathcal{U})$.

To complete our reduction of the theory of $\mathcal{H}$-modules for $\mathcal{M}_7(\mathcal{U})$ to that of modules for $\mathcal{U}$ we shall now show that every $\mathcal{H}$-module over $\mathcal{M}_7(\mathcal{U})$ is isomorphic to a module associated with a module for $\mathcal{U}$.

\begin{coro}\label{coro3} 
Let $V$ be a $\mathcal{H}$-module for $\mathcal{M}_7(\mathcal{U})$ such that $v\mathcal{M}_7(F)\neq 0$ for any $0\neq v\in V$. Then there exists a unital associative and commutative module $W$ for $\mathcal{U}$ such that $V$ is isomorphic to $\mathbf{T}(W)$.
\end{coro}
\begin{dem}
The fact that $V$ is a $\mathcal{H}$-module over $\mathcal{M}_7(\mathcal{U})$ such that $v\mathcal{M}_7(F)\neq 0$ for any $0\neq v\in V$,
implies that in particular $v\mathcal{M}_7(\mathcal{U})\neq 0$ for any $0\neq v\in V$. 
Let $\mathcal{M}=\mathcal{M}_7(\mathcal{U})\oplus V$ be
the split null extension of $\mathcal{M}_7(\mathcal{U})$ by $V$. So $\mathcal{M}$ is a Malcev algebra in $\mathcal{H}$ containing the 7-dimensional simple non-Lie Malcev algebra $\mathcal{M}_7(F)$ as a subalgebra such that $m\mathcal{M}_7(F)\neq 0$ for any $0\neq m\in\mathcal{M}$, then by Theorem $\ref{t4}$ there is a certain commutative and associative algebra $D$ such that 
$\mathcal{M}=\mathcal{M}_7(D).$
Let $W$ be the set of the elements of $D$ that appear in the matrix entries of $V$.
Thus
$V=\mathcal{M}_7(W)$
where $W\triangleleft D$ and $W^2=0$ in $D$, because $V\triangleleft \mathcal{M}$ and
$V^2=0$ in $\mathcal{M}$; so $D=\mathcal{U}\oplus W$. Then $D$ is a split null extension of $\mathcal{U}$ by its module $W$, hence $W$ is an associative and commutative module over $\mathcal{U}$ such that $w\mathcal{U}\neq 0$ for any $0\neq w\in W$. Therefore $\mathbf{T}(W)=V.$
\end{dem}

Evidently, $\mathbf{T}(\mbox{Reg}\,\mathcal{U})=
\mbox{Reg}(\mathcal{M}_7(\mathcal{U}))$. Is straightforward to proof that $\mathbf{T}$ is faithful and full, and using Corollary \ref{coro3} we have the desired equivalence of categories.

\begin{theo}\label{c3}
The categories $\textbf{\textup{Mod}}_\textbf{\textup{Assoc-Com}}\,\mathcal{U}$ and  $\mathcal{H}-\textbf{\textup{Mod}}_\textbf{\textup{Malc}}\,\mathcal{M}_7(\mathcal{U})$ are Morita equivalent.
\end{theo}

Now in the case of Malcev superalgebras, the Theorem \ref{t5} implies that if  $\mathbb{M}=\mathcal{M}_7(F)$, then the Malcev superalgebra $\mathcal{M}=\mathcal{M}_0\oplus \mathcal{M}_1$ in $\mathcal{H}$ is coordinated by $U=U_0\oplus U_1$, that is, $\mathcal{M}\cong\mathcal{M}_7 (U)$. Thus the following result establishes that in this case $\mathcal{M}=\mathcal{M}_0\oplus \mathcal{M}_1$ is Morita equivalent to $U=U_0\oplus U_1$. 

Let $U$ be a unital associative supercommutative superalgebra, and let $C(U)$ be the Cayley-Dickson matrix superalgebra over $U$. So introducing the supercommutator $[A,B]=AB-(-1)^{|A||B|}BA$ in $C(U)$ we obtain the Malcev superalgebra $C(U)^{(-)}$ and the non-Lie Malcev superalgebra
$\mathcal{M}_7(U)=C(U)^{(-)}/F\cdot e_{0}$.
Note that $\mathcal{M}_7(U)$ contains the 7-dimensional simple non-Lie Malcev algebra $\mathcal{M}_7(F)$.

Let's consider the following categories:
\begin{itemize}
\item let $\mathcal{H}-\widetilde{\textbf{\mbox{Mod}}}_\textbf{\mbox{Malc}}\,\mathcal{M}_7(U)$ denote the category of $\mathcal{H}$-modules $V=V_0\oplus V_1$ over $\mathcal{M}_7(U)$ such that $v\mathcal{M}_7(F)\neq 0$ for any homogeneous element $0\neq v\in V_0\cup V_1$.
\item let $\widetilde{\textbf{\mbox{Mod}}}_\textbf{\mbox{Assoc-Com}}U$ denote the category of unital associative supercommutative supermodules $W=W_0\oplus W_1$ over $U$ such that $wU\neq 0$ for any homogeneous element $0\neq w\in W_0\cup W_1$.
\end{itemize}
Thus analogous to the case of Malcev algebras we can define a functor 
$$\mathbf{S}:\widetilde{\textbf{\mbox{Mod}}}_\textbf{\mbox{Assoc-Com}}U\longrightarrow \mathcal{H}-\widetilde{\textbf{\mbox{Mod}}}_\textbf{\mbox{Malc}}\,\mathcal{M}_7(U)
.$$
We can get a result analogous to the Corollary \ref{coro3} and from there $\mathbf{S}(\mbox{Reg}\,U)=
\mbox{Reg}(\mathcal{M}_7(U))$. Also, is straightforward to proof that $\mathbf{S}$ is faithful and full. Then we have equivalence of categories.
 
\begin{theo}\label{c4}
The categories $\widetilde{\textbf{\textup{Mod}}}_\textbf{\textup{Assoc-Com}}U$ and   
$\mathcal{H-}\widetilde{\textbf{\textup{Mod}}}_\textbf{\textup{Malc}}\,\mathcal{M}_7(U)$
are Morita equivalent.
\end{theo}


\section{(Super)algebras and Modules with (super)involution}

Consider $\algA$ an arbitrary algebra. Recall that a linear mapping $\ast:\mathcal{A}\longrightarrow \mathcal{A}$ is called an \textit{involution} of an algebra $\mathcal{A}$, if it satisfies the conditions
$$(a^{\ast})^{\ast}=a,\,\, (ab)^{\ast}=b^{\ast}a^{\ast}$$
for any elements $a,b\in \mathcal{A}$.

If $\algA=\algA_0\oplus\algA_1$ is an arbitrary superalgebra. A linear even mapping $\ast:\mathcal{A}\longrightarrow \mathcal{A}$ is called a \textit{superinvolution} of a
superalgebra $\algA$, if it satisfies the conditions
$$(a^{\ast})^{\ast}=a,\,\, (ab)^{\ast}=(-1)^{\mid a\mid\mid b\mid}b^{\ast}a^{\ast}$$
for any homogeneous elements $a, b\in\algA_0\cup \algA_1$.

Now, let $V$ be a module over an algebra $(\mathcal{A},\ast)$ with involution. We will call $V$ an \textit{$\mathcal{A}$-module with involution}, if there exists a linear mapping
$-: V\longrightarrow V$ such that the mapping
$$v+a \mapsto \overline{v}+a^{\ast}$$
is an involution of the split null extension algebra $E=V+\mathcal{A}$. Evidently,
for an algebra with involution $\mathcal{A}$, the bimodules $\mbox{Reg}\, \mathcal{A}$ and $(\mbox{Reg}\, \mathcal{A})^{op}$ have the involutions induced by that of $\mathcal{A}$.

Let $\mathbb{O}$ be the Cayley-Dickson algebra over $F$. Then we know that $\mathbb{O}=F \oplus \mathbb{M}$, where 
$\mathbb{M}=\{x\in \mathbb{O} : t(x)=0 \}$ and multiplication in $\mathbb{O}$, for elements $a,b\in\mathbb{M}$ is defined by the following relation:
$$a\cdot b=-(a,b)+a\times b$$
where $(,)$ is a non-singular symmetric bilinear form on $\mathbb{M}$ and $\times$ is anticommutative multplication on $\mathbb{M}$. As $t(x)=x+\overline{x}$, then
$$\overline{a}=-a$$
for all $a\in\mathbb{M}$. So $\overline{\overline{a}}=a$ and 
$\overline{ab}=-ab\overset{(\ref{line1})}{=}ba=(-b)(-a)=\overline{b}\,\overline{a}$; hence $a\mapsto \overline{a}=-a$ is an involution in $\mathbb{M}$. 

So in this section we assume that in every Malcev (super)algebra that contains $\mathbb{M}$ is defined a non-singular symmetric bilinear form and an (super)involution that extend that of $\mathbb{M}$.

Let $V$ be a Malcev module for $\mathbb{M}$ and let $E=V+\mathbb{M}$ be the split null extension. So $E$ is a Malcev algebra that contains $\mathbb{M}$. Let 
$$f:E\times E\longrightarrow F$$ 
denote the non-singular symmetric bilinear form defined in $E$. 
Then consider the associative algebra with involution $\mathcal{L}(E)$, where
\begin{equation*}
\begin{split}
\mathcal{L}(E)&=\{ \alpha\in\mbox{End}_F(E): \,\exists\,\alpha^{*}\in\mbox{End}_F(E) \,\mbox{ such that}\\
&\,\,\,\,\,\,\,\,\,f(x\alpha,y)=f(x,y\alpha^{*})\,\,  \forall x,y\in E \}.
\end{split}
\end{equation*}
Note that $\alpha^{*}$ is unique by the non-singularity of $f$, which gives the involution in $\mathcal{L}(E)$.

Now we will study the structure of modules with involution over $\mathbb{M}$. Our objective is to prove that every module with involution over this algebra is completely reducible and every irreducible module with involution is of the type $\mbox{Reg}\, \mathbb{M}$.
In fact, we will consider modules with involution that satisfy the
additional condition of so-called \textit{J-admissibility} (see $\cite{J2}$). A module with
involution $(V,-)$ over an algebra with involution $(\mathcal{M}, \ast)$ is called \textit{J-admissible} if all the symmetric elements of the algebra with involution
$E = V+\mathcal{M}$ lie in the associative center (the nucleus) of E. 

\begin{theo}\label{t6}
Let $V$ be a $J$-admissible Malcev module over $\mathbb{M}$ and let $E=V\oplus \mathbb{M}$ be the corresponding split extension. 
If $f:E\times E\longrightarrow F$ is a non-singular symmetric bilinear form defined in $E$.
Then $V$ is completely reducible and
is a direct sum of irreducible modules with involution isomorphic to $\textup{Reg}\,\mathbb{M}$, that is, isomorphic to $\mathbb{M}\alpha$ for a certain symmetric element $\alpha$ of the centroid $\Gamma(E)$ with involution.
\end{theo}
\begin{dem}
Let $V$ be a module under consideration, with a involution $v\mapsto \overline{v}$. By Corollary $\ref{c1}$, $V$ is completely reducible, that is, $V=\sum_{i}\oplus V_{i}$, where $V_{i}$ is an irreducible almost faithful non-Lie Malcev $\mathbb{M}$-module. Also, by Theorem $\ref{t3}$ every irreducible almost faithful non-Lie Malcev $\mathbb{M}$-module is the regular module $\mbox{Reg}\,\mathbb{M}$. 
So we have $V_i=\mathbb{M}\alpha_{i}$ as modules, where $\alpha_{i}:\mathbb{M}\longrightarrow V_i$ is the isomorphism of modules that was proved in the Theorem $\ref{t3}$
given by $\alpha_i(z)=z\alpha_{i} (v_i,a_i,b_i)$ for some $v_i\in V_i$ and $a_i,b_i\in \mathbb{M}$ with $\alpha_{i} (v_i,a_i,b_i)\neq 0$. 
The condition $0\neq\alpha_{i} (v_i,a_i,b_i)$ is true since by Proposition $\ref{p3}$ 
$$0\neq\alpha_{i} (V_i,\mathbb{M},\mathbb{M})$$ because $\mathbb{M}$ is a non-Lie Malcev algebra.
So
$V=\sum_{i}\oplus \mathbb{M}\alpha_{i}$.

Denote by $\mathcal{U}=\sum_i F\alpha_{i}$ the span of all elements $\alpha_{i}$. 

By Theorem $\ref{t2}$ $E$ satisfies $h=0$, then by Lemma $\ref{l1}$ for all $x,y\in E$
\begin{equation*}
    p(x,v_{i},a_{i},b_{i})y=p(xy,v_{i},a_{i},b_{i}),
\end{equation*}
$(x\alpha_{i})y=(xy)\alpha_{i}$; thus $\alpha_{i}\in\Gamma(E)$, where $\Gamma(E)$ is \textit{the centroid} of $E$.

We can associate the involution $x\longmapsto \overline{x}$ of $E$ with the involution of $\mathcal{L}(E)$ by
\begin{equation}\label{e9}
\overline{x\alpha}=\overline{x}\alpha^{\ast}
\end{equation}
for all $x\in E$ and $\alpha\in \mathcal{L}(E)$. 

Fix arbitrary $a,b\in\mathbb{M}$, $v\in V$ and $\alpha\in \Gamma(E)$. Then there exist unique $\alpha^{\ast}\in\mathcal{L}(E)$ and so by $(\ref{e9})$
$$(ab)\alpha^{\ast}=(\overline{a}\,\overline{b})\alpha^{\ast}=\overline{(ba)}\alpha^{\ast}=\overline{(ba)\alpha}=\overline{(b\alpha)a}=
\overline{a}\overline{(b\alpha)}=a(b\alpha^{\ast}),$$
$$(ab)\alpha^{\ast}=(\overline{a}  \, \overline{b})\alpha^{\ast}=\overline{(ba)}\alpha^{\ast}=\overline{(ba)\alpha}=\overline{b(a\alpha)}=
\overline{(a\alpha)}\overline{b}=(a\alpha^{\ast})b,$$
$$(va)\alpha^{\ast}=-(\overline{a\overline{v}})\alpha^{\ast}=-\overline{(a\overline{v})\alpha}=-\overline{(a\alpha)\overline{v}}=-
\overline{\overline{v}}\overline{(a\alpha)}=v(x\alpha^{\ast})$$
and 
$$(va)\alpha^{\ast}=-(\overline{a\overline{v}})\alpha^{\ast}=-\overline{(a\overline{v})\alpha}=-\overline{a(\overline{v}\alpha)}=-
\overline{(\overline{v}\alpha)}\overline{a}=(v\alpha^{\ast})a;$$
thus $\Gamma(E)^{\ast}\subseteq \Gamma(E)$ and $\Gamma(E)^{\ast}= \Gamma(E)$. 
So $\Gamma(E)$ is a subalgebra with involution of $\mathcal{L}(E)$; hence  
$$\Gamma(E)=\textup{Sym}\,\Gamma(E)\oplus \textup{Skew}\,\Gamma(E).$$
Assume that there exists $0\neq \alpha\in \Gamma(E)$ such that $\alpha^{\ast}=-\alpha$. Then, if $\mathbb{M}=\mathcal{M}_7(F)$ we denote $a=e_1\alpha$, then $\overline{a}=a$; so $a$ is a symmetric element; hence, by $J$-admissibility of $V$ and using the mutiplication table of $\mathcal{M}_7(F)$ 
$$0=(e_2,e_1\alpha,e_3)=(e_2(e_1\alpha))e_3-e_2((e_1\alpha)e_3)=-8e_7\alpha$$ which implies $e_i\alpha=0$ for all $i=1,\dots,7$. Hence $\mathcal{M}_7(F)\alpha=0$
and $\alpha=0$, a contradiction. Similarly in the case $\mathbb{M}=\mathbb{M}_7$, the element $a=i\alpha$ is symmetric, then by $J$-admissibility of $V$ and using the mutiplication table of $\mathbb{M}_7$, we have $0=(j,i\alpha,k)=4(vk)\alpha$; so $\mathbb{M}_7\alpha=0$. Hence $\alpha=0$, a contradiction. Thus $\Gamma(E)=\textup{Sym}\,\Gamma(E).$

Now, if $V$ is irreducible then, let us take elements $v\in V$ and $a,b\in\mathbb{M}$ with $\alpha=\alpha(v,a,b)\neq 0$, where $\alpha\in\mathcal{U}$ and $V = \mathbb{M}\alpha$,
which is isomorphic to $\mbox{Reg}\, \mathbb{M}$, under
the isomorphism $\alpha:z\mapsto z\alpha(v,a,b)$.

In the general case, it suffices to notice that every $\alpha_i$ generates an irreducible module which is invariant under the involution and is isomorphic
to $\mbox{Reg}\, \mathbb{M}$.

The theorem is proved.
\end{dem}

The  following results affirm that the commutative associative algebra $\mathcal{U}$ of the theorem $\ref{t4}$ and the supercommutative associative superalgebra $U$ of the theorem $\ref{t5}$ are generated by symmetric elements. 

Let $\mathcal{M}$ be a Malcev algebra in $\mathcal{H}$ and let $f:\mathcal{M}\times \mathcal{M}\longrightarrow F$ denote the non-singular symmetric bilinear form defined in $\mathcal{M}$. As in the previous paragraphs, consider the associative algebra with involution $\mathcal{L}(\mathcal{M})$ with 
\begin{equation*}
\begin{split}
\mathcal{L}(\mathcal{M})&=\{ \alpha\in\mbox{End}_F(\mathcal{M}): \,\exists\,\alpha^{*}\in\mbox{End}_F(\mathcal{M}) \,\mbox{such that}\\
&\,\,\,\,\,\,\,\,\,f(x\alpha,y)=f(x,y\alpha^{*})\,\,  \forall x,y\in\mathcal{M}\}.
\end{split}
\end{equation*}
Also $\alpha^{*}$ is unique by the non-singularity of $f$, which gives the involution in $\mathcal{L}(\mathcal{M})$.

\begin{coro}\label{coro1}
Let $\mathcal{M}$ be a Malcev algebra in $\mathcal{H}$ with J-admissible involution (that is, every symmetric element lies in the nucleus of $\mathcal{M}$) such that $\mathcal{M}$ contains $\mathbb{M}$ as a subalgebra with $m\mathbb{M}\neq 0$ for any $m\neq 0$ from $\mathcal{M}$. 
If $f:\mathcal{M}\times \mathcal{M}\longrightarrow F$ is a non-singular symmetric bilinear form defined in $\mathcal{M}$.
Then $\mathcal{M}\cong\mathbb{M}\otimes_{F}\mathcal{U}$
for a certain commutative associative algebra $\mathcal{U}$ of symmetric elements.
\end{coro}
\begin{dem}
By Theorem \ref{t4} we have $\mathcal{M}\cong\mathbb{M}\otimes_{F} \mathcal{U}$, where $\mathcal{U}$ is a certain commutative subalgebra of $\Gamma(\mathcal{M})$. 
Using the same arguments from Theorem $\ref{t6}$, we obtain that
$\Gamma(\mathcal{M})$ is a subalgebra with involution of $(\mathcal{L}(\mathcal{M}),\ast)$ and is formed only by symmetric elements. Thus $\mathcal{U}$ is a subalgebra with involution of $\Gamma(\mathcal{M})$ whose elements are symmetric.

The corollary is proved.
\end{dem}

Now let $\mathcal{M}=\mathcal{M}_0\oplus\mathcal{M}_1$ be a Malcev superalgebra in $\mathcal{H}$
and let 
$h:\mathcal{M}\times \mathcal{M}\longrightarrow F$ denote the non-singular symmetric bilinear form defined in $\mathcal{M}$.
Consider the associative superalgebra with superinvolution $\mathfrak{L}(\mathcal{M})$ with
\begin{equation*}
\begin{split}
\mathfrak{L}(\mathcal{M})_i&=\{ \alpha\in\mbox{End}_F(\mathcal{M})_i: \,\exists\,\alpha^{*}\in\mbox{End}_F(\mathcal{M}) \,\mbox{ such that}\\
&\,\,\,\,\,\,\,\,\,h(x\alpha,y)=(-1)^{\mid \alpha\mid\mid y\mid}h(x,y\alpha^{*})\,\,  \forall x,y\in \mathcal{M}_0\cup\mathcal{M}_1 \}.
\end{split}
\end{equation*}
We observe that $\alpha^{*}$ is unique by the non-singularity of $h$, which gives the superinvolution in $\mathfrak{L}(\mathcal{M})$.

The following result is a generalization of the Corollary \ref{coro1} for superalgebras.

\begin{coro}\label{coro2}
Let $\mathcal{M}=\mathcal{M}_0\oplus\mathcal{M}_1$ be a Malcev superalgebra in $\mathcal{H}$ with J-admissible superinvolution (that is, every symmetric element lies in the nucleus of $\mathcal{M}$) such that $\mathcal{M}_0$ contains $\mathbb{M}$, with $m\mathbb{M}\neq 0$ for any $m\neq 0$ from $\mathcal{M}_0\cup\mathcal{M}_1$. 
If $h:\mathcal{M}\times \mathcal{M}\longrightarrow F$ is a non-singular symmetric bilinear form defined in $\mathcal{M}$.
Then $\mathcal{M}\cong\mathbb{M}\otimes_{F} U$
for a certain supercommutative and associative superalgebra $U$ of symmetric elements.
\end{coro}
\begin{dem}
By Theorem \ref{t5} we have $\mathcal{M}\cong\mathbb{M}\otimes_F U$, where $U=U_0+ U_1$ is a certain supercommutative subsuperalgebra of $\widetilde{\Gamma}(\mathcal{M})$. 

Consider the superalgebra $\mathfrak{L}(\mathcal{M})$ with superinvolution $\ast$. We can associate the superinvolution $a\mapsto \overline{a}$ of $\mathcal{M}$ with the superinvolution of $\mathfrak{L}(\mathcal{M})$ by
\begin{equation}\label{e20}
\overline{a\alpha}=\overline{a}\alpha^{\ast}
\end{equation}
for all $a\in\mathcal{M}$ and $\alpha\in \mathfrak{L}(\mathcal{M})$. 

Remember that as $\widetilde{\Gamma}(\mathcal{M})$ is a subsuperalgebra of $\mbox{End}_F(\mathcal{M})$ we have 
$\mathcal{M}_i\alpha_j\subseteq\mathcal{M}_{i+j}$ for any homogeneous elements $\alpha_j$ of $\widetilde{\Gamma}(\mathcal{M})$. Fix arbitrary $x,y\in\mathcal{M}$ and $\alpha\in\widetilde{\Gamma}(\mathcal{M})$. Then there exist unique $\alpha^{\ast}\in\mathfrak{L}(\mathcal{M})$ such that
$$(xy)\alpha^{\ast}=(\overline{\overline{x}}\,\overline{\overline{y}})\alpha^{\ast}=
(-1)^{\mid x\mid\mid y\mid}\overline{(\overline{y}\,\overline{x})}\alpha^{\ast}
\overset{(\ref{e20})}{=}(-1)^{\mid x\mid\mid y\mid}\overline{(\overline{y}\,\overline{x})\alpha}:$$
\begin{equation*}
\begin{split}
(-1)^{\mid x\mid\mid y\mid}\overline{(\overline{y}\,\overline{x})\alpha}&=(-1)^{\mid x\mid\mid y\mid}(-1)^{\mid x\mid\mid\alpha\mid}
\overline{(\overline{y}\alpha)\overline{x}},\\
&=(-1)^{\mid x\mid(\mid y\mid+\mid\alpha\mid)}(-1)^{\mid x\mid\mid \overline{y}\alpha\mid}\overline{\overline{x}}\,\overline{\overline{y}\alpha}\overset{(\ref{e20})}{=}(-1)^{\mid x\mid(\mid y\mid+\mid\alpha\mid)}(-1)^{\mid x\mid(\mid y\mid+\mid\alpha\mid)}x(y\alpha^{\ast}),\\
&=x(y\alpha^{\ast})
\end{split}
\end{equation*}
and
\begin{equation*}
\begin{split}
(-1)^{\mid x\mid\mid y\mid}\overline{(\overline{y}\,\overline{x})\alpha}&=(-1)^{\mid x\mid\mid y\mid}
\overline{\overline{y}(\overline{x}\alpha)},\\
&=(-1)^{\mid x\mid\mid y\mid}(-1)^{\mid \overline{x}\alpha\mid\mid y\mid}\overline{\overline{x}\alpha}\,\overline{\overline{y}}
\overset{(\ref{e20})}{=}(-1)^{\mid x\mid\mid y\mid}(-1)^{(\mid x\mid+\mid\alpha\mid)\mid y\mid}(x\alpha^{\ast})y,\\
&=(-1)^{\mid\alpha\mid\mid y\mid}(x\alpha^{\ast})y.
\end{split}
\end{equation*}
Then
$(xy)\alpha^{\ast}=x(y\alpha^{\ast})=(-1)^{\mid\alpha\mid\mid y\mid}(x\alpha^{\ast})y;$ 
thus $\widetilde{\Gamma}(\mathcal{M})^{\ast}\subseteq \widetilde{\Gamma}(\mathcal{M})$ and $\widetilde{\Gamma}(\mathcal{M})^{\ast}= \widetilde{\Gamma}(\mathcal{M})$. 
So $\widetilde{\Gamma}(\mathcal{M})$ is a subsuperalgebra with superinvolution of $\mathfrak{L}(\mathcal{M})$ and hence it is clear that
$$\widetilde{\Gamma}(\mathcal{M})=\textup{Sym}\,\widetilde{\Gamma}(\mathcal{M})\oplus \textup{Skew}\,\widetilde{\Gamma}(\mathcal{M}).$$
Assume that there exists $0\neq \alpha\in \widetilde{\Gamma}(\mathcal{M})$ such that $\alpha^{\ast}=-\alpha$. So as in the Theorem $\ref{t6}$ if $\mathbb{M}=\mathcal{M}_7(F)$ we denote $a=e_2\alpha$ and $\overline{a}=a$; then $a$ is a symmetric element and by $J$-admissibility of $\mathcal{M}$ and using the mutiplication table of $\mathcal{M}_7(F)$ we have
$0=(e_4,e_2\alpha,e_3)=(e_4(e_2\alpha))e_3-e_4((e_2\alpha)e_3)=-4e_1\alpha$ which implies $e_i\alpha=0$ for all $i=1,\dots,7$. Hence $\mathcal{M}_7(F)\alpha=0$
and $\alpha=0$, a contradiction. Similarly in the case $\mathbb{M}=\mathbb{M}_7$ we get a contradiction. Therefore $\widetilde{\Gamma}(\mathcal{M})=\textup{Sym}\,\widetilde{\Gamma}(\mathcal{M}).$

It is easy to see that $U$ is invariant under the superinvolution of $\widetilde{\Gamma}(\mathcal{M})$; thus $U$ is generated by symmetric elements.  

The corollary is proved.
\end{dem}

\section{Acknowledgement}


This work was funded by CONCYTEC-FONDECYT within the framework of the call “Proyecto Investigación Básica 2019-01” [380-2019-FONDECYT].

\section{Declarations}

\textbf{Conflict of Interests} The author declares that they have no conflict of interest.

\textbf{Availability of data and materials} The author declares that data supporting the findings of this study are available within the article.


\begin{thebibliography}{99}

\frenchspacing

\bibitem{E} {\em Elduque, A.}, On a class of Malcev superalgebras. J. Algebra 173 (1995) no. 2, 237–252.

\bibitem{ES1} {\em A. Elduque, I. P. Shestakov,} On Malcev superalgebras with trivial Lie nucleus, Nova J. Algebra Geom., 2 (1993), 361-366.

\bibitem{ES} {\em Elduque, A., Shestakov, I.P.}, Irreducible non-Lie modules for Malcev superalgebras, J. Algebra 173 (1995) 622–637

\bibitem{F2} {\em V. T. Filippov,} Nilpotent ideals of Mal'tsev algebras, Algebra and Logic 18 (1979), 379-389.

\bibitem{F4} {\em V. T. Filippov,} Finitely generated Mal'tsev algebras, Algebra and Logic 19 (1980), 480-499.

\bibitem{F} {\em V. T, Filippov}, Varieties of Mal'tsev algebras. Algebra and Logic 20 (1981) no 3, 200-210.

\bibitem{F1} {\em V. T. Filippov}, Prime Malcev algebras, Mat. Zametki 31 (1982), 669-677; English transl. in Math.
Notes 31 (1982).

\bibitem{F3} {\em V. T. Filippov,} Imbedding of Mal'tsev algebras into alternative algebras, Algebra and Logic 22 (1983), 443-465.

\bibitem{J1} {\em N.\,Jacobson}, A Kronecker factorization theorem for Cayley algebras and the exceptional simple Jordan algebra. Amer. J. Math. 76, (1954). 447-452.

\bibitem{J2} {\em N. Jacobson}, Structure and Representations of Jordan Algebras, Amer. Math. Soc. Colloq. Publ., Vol. XXXIX, Amer. Math. Soc., Providence, RI, 1968

\bibitem{K} {\em I.\,Kaplansky}, Semi-simple alternative rings, Portugal. Math.10 (1966), 37-50.

\bibitem{LDS1} {\em López-Díaz,\,M.\,C., Shestakov, Ivan P.}, Representations of exceptional simple Jordan superalgebras of characteristic 3. Comm. Algebra 33 (2005), no. 1, 331-337.

\bibitem{LDS2} {\em López-Díaz,\,M.\,C., Shestakov, Ivan P.}, Representations of exceptional simple alternative superalgebras of characteristic 3. Trans. Amer. Math. Soc. 354 (2002), no. 7, 2745-2758.

\bibitem{LSS} {\em López-Solís,\,V.\,H.}, Kronecker Factorization Theorems for Alternative Superalgebras, J. Algebra 528 (2019) 311-338.

\bibitem{LSS1} {\em Solís, Victor Hugo López, and Ivan P. Shestakov.}, On a problem by Nathan Jacobson. Revista Matemática Iberoamericana (2021).

\bibitem{CE} {\em Martínez,\,C.,  Zelmanov, E.}, A Kronecker factorization theorem for the exceptional Jordan superalgebra. Journal of Pure and Applied Algebra 177 (2003) no.1, 71-78.

\bibitem{PSS} {Pchelintsev, \,S. \,V., Shashkov, \,O. \,V. and Shestakov, \,I. \,P.}, Right alternative bimodules over Cayley algebra and coordinatization theorem, J. of Algebra 572 (2021) 111-128.

\bibitem{PY} {\em Popov, Y.}, Representations of simple noncommutative Jordan superalgebras I. J. Algebra 544 (2020): 329-390.

\bibitem{PS} {Pozhidaev, A. P., Shestakov, I. P. (2009).} Noncommutative Jordan superalgebras of degree n > 2, Algebra Log. 49 (1) (2010) 18–42.

\bibitem{SA} {\em Sagle, Arthur A.}, Malcev algebras. Trans. Amer. Math. Soc. 101 (1961) no 3, p. 426-458.

\bibitem{Sh} {\em Shestakov, I.P.}, Prime Malcev Superalgebras, Mat. Sb. 182 (1991), 1357–1366.


\end{thebibliography}
\end{document}